\documentclass[twocolumn]{svjour3}

\usepackage[english]{babel}     
\usepackage{lmodern}            
\usepackage[T1]{fontenc}        
\usepackage[utf8]{inputenc}     
\usepackage{amsmath}            
\usepackage{graphicx}
\usepackage{epstopdf}           
\usepackage{flushend}           
\usepackage[hidelinks, bookmarksdepth=3]{hyperref}           
\usepackage[round]{natbib}
\usepackage{indentfirst}
\usepackage{multirow}
\usepackage[caption=false]{subfig}
\usepackage{tikz}
\usepackage{pgfplots}
\pgfplotsset{compat=1.14}
\usepackage{booktabs}
\usepackage{tabularx}
\usepackage{array}
\usepackage{xspace}
\usepackage[title]{appendix}

\usepackage{amssymb}

\usepackage[figuresright]{rotating}

\usepackage{lipsum}
\usepackage{tabularx}
\usepackage[restart]{parnotes}

\newcommand{\xis}{\mathbf{x}}
\newcommand{\Xis}{\mathbf{X}}
\newcommand{\de}{\mathbf{d}}
\newcommand{\pdf}[1]{\mathit{f}_#1}

\newcommand{\deq}{\mathrel{\mathop:}=}
\newcommand{\ie}{\emph{i.e.}\@\xspace}
\newcommand{\eg}{\emph{e.g.}\@\xspace}
\newcommand{\etc}{\emph{etc.}\@\xspace}
\newcommand{\argmin}{\hspace{2pt}\textrm{argmin}}
\newcommand{\Cov}{\hspace{2pt}\textrm{Cov}}
\newcommand{\E}{\hspace{2pt}\textrm{E}}

\newcolumntype{C}{ >{\centering\arraybackslash} m{4cm} }
\newcolumntype{D}{ >{\centering\arraybackslash} m{1.3cm} }
\def\block(#1,#2)#3{\multicolumn{#2}{c}{\multirow{#1}{*}{$ #3 $}}} 

\begin{document}

\def\figureautorefname{Fig.}
\def\equationautorefname~#1\null{Eq.~(#1)\null}

\title{Monte Carlo Integration with adaptive variance selection for improved stochastic Efficient Global Optimization}

\author{Felipe Carraro \and Rafael Holdorf Lopez \and Leandro Fleck Fadel Miguel \and Andr\'e Jacomel Torii}

\institute{Felipe Carraro, \email{felipecarraro@gmail.com}\\ Rafael Holdorf Lopez,
\email{rafael.holdorf@ufsc.br}\\ Leandro Fleck Fadel Miguel, \email{leandro.miguel@ufsc.br}\\ \at
Center for Optimization and Reliability in Engineering (CORE), Civil Engineering Department, Federal
University of Santa Catarina, Rua João Pio Duarte Silva, s/n 88040-900 Florianópolis, SC, Brazil.
\and Andr\'e Jacomel Torii, \email{ajtorii@hotmail.com} \\\at
Center for Optimization and Reliability in Engineering (CORE), Civil Engineering Department, Federal University for Latin American Integration,  Av. Silvio Américo Sasdelli, 1842 - Vila A, 85866-000 Foz do Iguaçu, PR, Brazil
 }


\maketitle

\begin{abstract}
In this paper, the minimization of computational cost on evaluating multi-dimensional integrals is explored. More specifically, a method based on an adaptive scheme for error variance selection in Monte Carlo integration (MCI) is presented. It uses a stochastic Efficient Global Optimization (sEGO) framework to guide the optimization search. The MCI is employed to approximate the integrals, because it provides the variance of the error in the integration. In the proposed approach, the variance of the integration error is included into a Stochastic Kriging framework by setting a target variance in the MCI. We show that the variance of the error of the MCI may be controlled by the designer and that its value strongly influences the computational cost and the exploration ability of the optimization process. Hence, we propose an adaptive scheme for automatic selection of the target variance during the sEGO search. The robustness and efficiency of the proposed adaptive approach were evaluated on global optimization stochastic benchmark functions as well as on a tuned mass damper design problem. The results showed that the proposed adaptive approach consistently outperformed the constant approach and a multi-start optimization method. Moreover, the use of MCI enabled the method application in problems with high number of stochastic dimensions. On the other hand, the main limitation of the method is inherited from sEGO coupled with the Kriging metamodel: the efficiency of the approach is reduced when the number of design variables increases.
\end{abstract}


\section{Introduction}

Much of today’s engineering analysis consists of running complex computer codes.
Although computer power has been steadily increasing, the expense of running many analysis codes
remain non-trivial. For example, a single evaluation of a finite element model may
take minutes to
hours, if not longer. Thus, it often becomes impractical to perform a large number of such
simulations, for example, as required during optimization. To address such a challenge,
approximation or metamodeling techniques are often used.

One of such techniques is Kriging \citep{Krige51,cressie1993statistics}, which has
been gaining popularity in the last decade, and has been applied in numerous branches of science
\citep{theodossiou2006evaluation,Goel2008,huang2011optimal,Sakata2011,zhu2014geometrically}.
It has been shown to produce accurate response surfaces and to offer a estimate of the error in points that had not been sampled.

Kriging in its usual formulation considers the approximation of deterministic
functions. But what if the objective function possess some sort of randomness or variability? This
randomness could be, for example, uncertainty in the input parameters or noise
in the function response. For such cases, a more recent extension called
Stochastic Kriging (SK)  \citep{staum2009better,ankenman2010stochastic,kleijnen2016estimating,kaminski2015method,chen2014stochastic, plumlee2014building} was developed. Improvements and extensions emerged for both SK and deterministic Kriging in recent years. For example, \citet{wang2014new} included qualitative factors on the SK metamodel by constructing correlation functions that are valid across levels of qualitative factors. 
Several authors incorporated gradient estimators into the metamodel, which tends to significantly improve surface prediction \citep{chen2013enhancing,qu2014gradient,ulaganathan2014use,hao2018adaptive}. \citet{chen2014stochastic} evaluated sampling strategies while taking into account possible bias in SK simulation response
estimates. \citet{plumlee2014building} suggested an approach where
the intrinsic metamodel of SK error is not assumed Gaussian distributed, deriving a distributional emulator called Quantile Kriging. 
In the same subject, \citet{zhang2017asymmetric} built the Asymmetric Kriging emulator, again without error distribution assumption. This approach employed 
asymmetric least squares, which are used to compute multiple quantile curves and fit the model. \citet{kaminski2015method} proposed an improved SK update procedure by combining metamodel and simulation evaluations in addition to the use of smoothed variance evaluations to
compute the intrinsic noise. \citet{shen2018enhancing} applied stylized queuing models to provide information about the shape of the response surface.
Following \citet{chen2012effects}, \citet{2018arXiv180200677Z} further explored the use of Common Random Numbers with SK, showing that its use may not always be detrimental to prediction accuracy, being beneficial if the observation errors of the real system are small enough.

A powerful approach to use Kriging or SK within an optimization context is the Efficient Global Optimization (EGO). EGO is the name of the framework developed by \citet
{Jones1998}, which exploits the information provided by the Kriging metamodel to iteratively add new
points, improving the surrogate accuracy and at the same time seeking its global minimum. It is the manner in which these infill points are added that characterizes different EGO methods.

The use of SK within the EGO framework, which we name here stochastic Efficient Global Optimization (sEGO), is relatively recent. For example, 
\citet{JALALI2017279} compared Kriging-based methods in heterogeneous noise situations, while \citet{picheny2013benchmark} benchmarked different infill criteria for the noisy case. The choice of infill in stochastic envioroment has recieved considerable attention in literature. Some authors considered the use of the classical (noiseless) 
Expected Improvement (EI) criterion with a ``plug-in'', which is effectively an target for the definition of improvement, such as \citet{vazquez2008global} and \citet{osborne2009gaussian}. Unfortunately, for the noisy case plug-ins do not take into acocunt the noise of future observations. The Augmented Expected Improvement (AEI) method from \citet{huang2006global}, adresses this issue with a penalty term. Such formulation is employed in this paper and is further discussed in \autoref{sub:augmented_expected_improvement}. \citet{forrester2008engineering} employed a reinterpolation procedure by building one model for evaluating and predicting noisy observations and another interpolative model for determining the infill point. \citet{picheny2013quantile} considered a modificition of EI using quantiles to define a reference for improvement, taking into account the effect of a new observation on the model. The work of \citet{scott2011correlated} applied the concept of knowledge gradient to define improvement, aiming to measure the global effect of a new measurement on the model mean. Methods employing a global measure of improvement have also been proposed, however those involve potentially expensive numerical integrations. \citet{villemonteix2009informational} criterion maximized information gain, while \citet{gramacy} evaluated the effect of a candidate point over the design space with an Integrad Expected Conditional Improvement.

In this context, this paper presents an efficient sEGO approach for optimization problems whose objective function depends on an integral. The proposed approach is based on Monte Carlo Integration (MCI) and sEGO. First, MCI is employed to approximate the objective function since it provides not only an approximation for the integral, but also the variance of the error. The variance of the error is then included into the SK framework, and the AEI infill criterion is employed to guide the addition of new points in the stochastic EGO framework. We show that the variance of the error committed by MCI may be controlled by the designer and that its value strongly influences the computational cost and the exploration ability of the optimization process. Thus, the main contribution of this paper is the development of an adaptive scheme for the evaluation of a target variance of the error committed by MCI. The adaptive approach provides a framework to avoid getting stuck in local minima due to lack of information, and to achieve an efficient optimization process by rationally spending the available computational budget, \ie it balances the exploration and exploitation capabilities of the algorithm.

The rest of the paper is organized as follows: \autoref{ModelSituation} presents the problem statement, \ie the characteristics of the problems we aim to solve as well as a description of the MCI. The sEGO framework is presented in \autoref{sec:ego}. The potential consequences of the target variance setting are detailed in \autoref{sec:howtosetvariance}, while the adaptive target approach is presented in \autoref{sec:adap_target}. Numerical examples are studied in \autoref{sec:numerical} to show the efficiency and robustness of the proposed method. Finally, the main conclusion are listed in \autoref{sec:conclusions}.

\section{Problem Statement}\label{ModelSituation}

The goal of this paper is to solve the problem of minimization of a function
$J$, which depends on an integral as in

\begin{equation}
 \underset{d \in S}{\min J(\de)} =  \int_{\Omega} \phi(\de,\textbf{x}) w(\mathbf{x}) d\textbf{x},
 \label{eq:problemabasico}
\end{equation}

\noindent where $\de \in \Re^n$ is the design vector, $\xis  \in \Re^{n_x}$ is the parameter vector, $\phi : \Re^n \times \Re^{n_x} \to \Re$ is a known function, $S$ is the design domain, $w(\mathbf{x})$ is some known weight function (\eg probability distribution) and $\Omega \subseteq \Re^{n_x}$ is the integration domain (\eg support of the probability distribution). We also assume here that the design domain $S$ considers only box constrains.

Here we are interested in situations that: 

\begin{itemize}
    \item this integral cannot be evaluated analytically,
    \item $\phi$ is a black box function and is computationally demanding, 
    \item the resulting objective function $J$ is not convex and multimodal,
    \item $S$ includes only box constrains.
\end{itemize}
 
To give only a few examples, in engineering applications problems of this kind may arise from the maximization of the expected performance of a mechanical system, vastly applied in robust design \citep{Capiez08a,Capiez08b,Ritto11,Lopez2014,leticia2016,MIGUEL2016703}, the multidimensional integral of Performance Based Design Optimization \citep{beck2014optimal,spence2014performance,Bobby2016}, or the double integral of Optimal Design of Experiment problems \citep{huan2013simulation,Beck2017fast}.

In the formulation of \autoref{eq:problemabasico}, it is important to notice the difference between $n$, which is the dimension of the optimization problem and $n_x$, which is the stochastic dimension of the integral. Indeed, in cases where the integration domain $\Omega$ of \autoref{eq:problemabasico} has high dimension, evaluation of the above integral may become computationally demanding. Numerical integration such as quadrature procedures may be employed, although its computational efficiency is drastically reduced for high dimensional problems. In such cases, together with the sEGO approach proposed in this paper, sampling techniques may be employed, such as MCI, Importance Sampling \citep{Rubinstein2007}, Multi Level Monte Carlo \citep{giles2008multilevel}, Multi index Monte Carlo \citep{haji2016multi}.

In MCI, $J$ is estimated as

\begin{equation}
    J(\de) \approx \bar{J}(\de) = \frac{1}{n_r} \sum_{i=1}^{n_r} \phi
    (\de,
    \xis^{(i)} ),
    \label{eq:approxmean}
\end{equation}

\noindent where $n_r$ is the sample size and $\xis^{(i)}$ are sample points randomly drawn from distribution $w(\mathbf{x})$. In this work the sample points employed for MCI are also called replications, in order to distinguish from sample points employed for Kriging, when necessary.

One of the advantages of using sampling techniques, such as MCI, is that we are able to estimate the variance of the error of the approximation. Indeed, the variance of the estimator can be computed for a fixed $\de$ by a point estimate as:

\begin{equation}
    \overline{\sigma}^{2}(\de) = \frac{1}{n_r(n_r-1)} \sum_{i=1}^{n_r} 
    (\phi_i -
    \bar{J}(\de))^2,
    \label{eq:pointestimate}
\end{equation}

\noindent where $\phi_i = \phi(\de,\xis^{(i)})$. Thus, by increasing the sample size $n_r$ (\ie the number of replications), the variance estimate decreases and approximation in \autoref{eq:approxmean} gets closer to the exact value of
\autoref{eq:problemabasico}. In fact, it may be demonstrated that the rate of convergence of MCI is $n_{r}^{-\frac{1}{2}}$ \citep{hammersley1964general,kalos1986monte}.

As already mentioned, in this paper we aim to solve problems in which $\phi$ is a black box computationally demanding function, and $J$ is nonconvex. We assume
gradient information is unknown, although such information, if available, could be employed as in \citet{hao2018adaptive}.
 EGO algorithms are able to handle these difficulties and were successfully applied in different fields \citep{couckuyt2010surrogate,li2011review,bae2012ego,duvigneau2012kriging,gengembre2012kriging,kanazaki2015efficient,Chaudhuri2015,Haftka2016,UrRehman2017}. More specifically, in this paper we propose the use of sEGO, but the SK metamodel is enriched with the information given by \autoref{eq:pointestimate}. EGO, SK and the proposed scheme are detailed in the next sections.

\section{Stochastic Efficient Global Optimization (sEGO)}
\label{sec:ego}

 According to \citet{Jones2001}, EGO methods generally follow these steps:

\begin{enumerate}
    \item Construction of the initial sampling plan;
    \item Construction of the metamodel;
    \item Addition of a new infill point to the sampling plan and return to step 2.
\end{enumerate}

Steps 2 and 3 are repeated until a stop criterion is met, \eg, maximum
number of function evaluations. The manner in which the infill points are added in each iteration is what differs the different EGO approaches. In the next subsections, these steps are detailed in order to set the basis of the proposed approach.

\subsection{Initial sampling plan}


In the first step, a Kriging sampling plan $\boldsymbol{\Gamma}$ containing $n_s$ points is created, \ie

\begin{equation}
\boldsymbol{\Gamma} = \{\de^{(1)},\de^{(2)},\dots,\de^{(n_s)}\}.
\end{equation}

\noindent A Latin Hypercube scheme is usually employed for this purpose. Then, the objective function value $J$ of each of these points is evaluated  using the original model, obtaining

\begin{equation}
\textbf{y} = \{y^{(1)},y^{(2)},\dots,y^{(n_s)}\},
\end{equation}

\noindent where $y^{(i)} = J\left(\de^{(i)}\right)$. 

Step 2 constructs a prediction model, which is given by the SK in this paper. In order to keep the paper self-contained, the formulation of Deterministic Kriging is given in Appendix \autoref{dKriging}, while the SK formulation is given in the next subsection.

\subsection{Stochastic Kriging (SK)}
\label{sec:stochastic_kriging}

\citet{ankenman2010stochastic} proposed an extension to the deterministic
Kriging methodology to deal with stochastic simulation. Their main
contribution was to  account for the sampling variability that is
inherent to a stochastic simulation.  In order to accomplish this, they
characterized both the intrinsic error inherent in a stochastic simulation
and the extrinsic error that comes from the metamodel approximation. Then, the SK prediction can be seen as:

\begin{align}
    \hat{y}(\de_i) =\overbrace{M(\de_i)}^{\text{Trend}} + \overbrace{Z(\de_i)}^{\text{Extrinsic}} + \overbrace{\epsilon(\de_i)}^{\text{Intrinsic}}, \nonumber\\ \quad 
\end{align}

\noindent where $M(\de)$ is the usual average trend, $Z(\de)$ accounts for the model uncertainty and
is now referred as extrinsic noise. The additional term $\epsilon$, called intrinsic noise,
accounts for the simulation uncertainty or variability. In the original version of \citet{ankenman2010stochastic}, the intrinsic noise is assumed independent and identically distributed (i.i.d.) across Kriging sample points and possess a Gaussian distribution with zero mean. 

In the present paper, the variability is due to the error in the approximation of the integral from Eq. (\ref{eq:problemabasico}) caused by MCI. It is worth to recall here that MCI provides an estimation of the variance of this error. That is, we are able to estimate the intrinsic noise, and consequently, introduce this information into the metamodel framework. In order to accomplish this,  we construct the covariance matrix of the intrinsic noise - among the current sampling plan points. Since the intrinsic error is assumed to be i.i.d. and Normal, the covariance matrix is a diagonal matrix with components

\begin{equation}
    (\Sigma_\epsilon)_{ii} = \overline{\sigma}^{2}(\de_i), 
    \quad i=1,2,...,n_s,
    \label{eq:diagonal_matrix}
\end{equation}

\noindent where $\overline{\sigma}^{2}$ is given by \autoref{eq:pointestimate}. Then, considering the Best Linear
Unbiased Predictor shown by \citet {ankenman2010stochastic}, the prediction of the SK  at a given point $\de_u$ is:

\begin{equation}
    \hat{y}(\de_u) =  \hat{\mu} + \textbf{r}^T (\boldsymbol{\Psi}+\Sigma_\epsilon)^{-1} (\textbf{y} -
    \textbf{1}
    \hat{\mu} ),
    \label{predSK}
\end{equation}

\noindent which is the usual Kriging prediction with the added
diagonal correlation matrix from the intrinsic noise, as can be seen by comparing Eqs.  (\ref{predSK}) and (\ref{eq:krigingpredictor}) .

Similarly, the predicted error takes the form:


\begin{equation}
\begin{aligned}
s_n^2(\de)  = & \widehat{\sigma^2} \left[ 1 + \lambda(\de) - \textbf{r}^T (\boldsymbol{\Psi}+\Sigma_\epsilon)^{-1} \textbf{r} \right.\\
              & + \left. \frac{(1 - \textbf{1}^T(\boldsymbol{\Psi}+\Sigma_\epsilon)^{-1} \textbf{r})^2}{\textbf{1}^T (\boldsymbol{\Psi}+\Sigma_\epsilon)^{-1} \textbf{1}} \right],
\end{aligned}
\label{eq:skestimatederror}
\end{equation}

\noindent where $\lambda(\de)$ corresponds to the regression term. If
the errors are assumed homoscedastic, then this term becomes a
constant and $\Sigma_\epsilon = I \lambda$. In the present paper, the designer may set its value in the MCI, \ie its value depend on the selected variance 
 $\overline{\sigma}^{2}$. Further discussion on the how targets are chosen is
presented in \autoref{sec:howtosetvariance}.

\subsection{Infill criterion for stochastic EGO}
\label{sec:infill_criteria_for_noisy_evaluations}
\label{sub:augmented_expected_improvement}

Although SK provides a good noise filtering model, the error estimates are no longer appropriate for use when choosing
infill points in the EGO framework \citep{j2006design}. 

The popular Expected Improvement (EI) infill criteria has proven
 global convergence under the deterministic case \citep{Locatelli1997}. However, this
 property no longer holds for the stochastic case, which does not
 guarantee a dense sampling. Overall, any method that relies on the
 predicted error going to zero at sampled points, which is the case in the
 deterministic approach, has difficulties for application in
 stochastic environments \citep{j2006design}. Hence, an infill point method adapted to SK is required. \citet{picheny2013benchmark} benchmarked different infill criteria for SK case. From that paper, a modification of
the deterministic EI criterion called AEI \citep{huang2006global} provided promising results and it is employed here. The main steps of the AEI criterion are presented in the next paragraphs, while a  full description is given by \citet{huang2006global} and \citet{picheny2013benchmark}. Moreover, comparison of AEI to deterministic infill criteria may be find in references \citet{huang2006global,forrester2008engineering}.

In the EI criterion, the infill point is selected as the one that maximizes the
expected value of the improvement measurement. This improvement definition requires a target value
to indicate the greediness of the search. Considering the improvement definition:

\begin{equation}
    I(\de, y_{target}) \deq \max(0,y_{target} - Y(\de)),
\end{equation}

\noindent where the $y_{target}$ for EI is usually chosen as the minimal solution found so far. For the AEI
infill criterion, this target is the so-called effective best solution $\de^{**}$ and is computed as:

\begin{equation}
    \de^{**} = \argmin(\de^{(i)} + \alpha s^{(i)}_n) \quad \text{for }i=1,...,n_s,
\end{equation}\\

\noindent  where  $s^{(i)}_n$ is the corresponding
kriging error, obtained by taking the square root of the MSE defined in \autoref{eq:skestimatederror}, and
$\alpha$ is an arbitrary constant. 

The criterion can be
calculated as the expected improvement over the effective best solution multiplied by a penalization
term:

\begin{equation}
    AEI(\de) = \E[I(\de, \de^{**})] \left( 1 - \frac{\sqrt{\lambda(\de)}}{\sqrt{s_n^2(\de) + \lambda(\de)}} \right),
    \label{AEIinfill}
\end{equation}

\noindent where $\E[I(\de, \de^{**})]$ represents the expected value of
the already defined improvement, $s_n^2$ the Kriging estimated error
and $\lambda$ the intrinsic output noise, which will be enforced by
MCI procedure based on a target variance.


The term between parenthesis in the right-hand-side of Eq. (\ref{AEIinfill}) is a penalty factor which amplifies the importance of the Kriging variance. That is, it enhances exploration, avoiding multiple simulations over the same input \citep{picheny2014noisy}. 


\section{Problem: how to set the target variance?}
\label{sec:howtosetvariance}

With the framework presented so far, we are able to incorporate error estimates from MCI within the sEGO scheme. It is important to notice that the number of samples (replications) of the MCI is an input parameter, \ie the designer has to set $n_r$ in Eq. (\ref{eq:pointestimate}). Consequently, the designer is able to control the magnitude of $\Sigma_\epsilon$ and $\lambda$ by changing the sample size $n_r$. However, in practice a target variance ($\overline{\sigma}^{2}_{target}$) is first chosen and the sample size is iteratively increased until the evaluated variance is close to the target value. Thus, for a constant target variance, the regression parameter is then enforced by the MCI procedure to be

\begin{equation}
    \lambda(\de) = \overline{\sigma}^{2}_{target}.
\end{equation}

The choice of the target variance must consider two facts:

\begin{itemize}
\item if the target variance is too high, the associated error may lead to a poor and deceiving approximation of the integral,
\item if the target tends to zero, so does the error and we retrieve the deterministic case, however, at the expense of a huge computational effort.
\end{itemize}
 
This section aims to emphasize the importance of the target variance setting and its consequences in the optimization process using stochastic EGO. It presents the key ideas behind the main contribution of this paper: the adaptive target variance approach proposed in Section \ref{sec:adap_target}.

\begin{figure}[htbp]
    \centering
    \includegraphics[width=\linewidth]{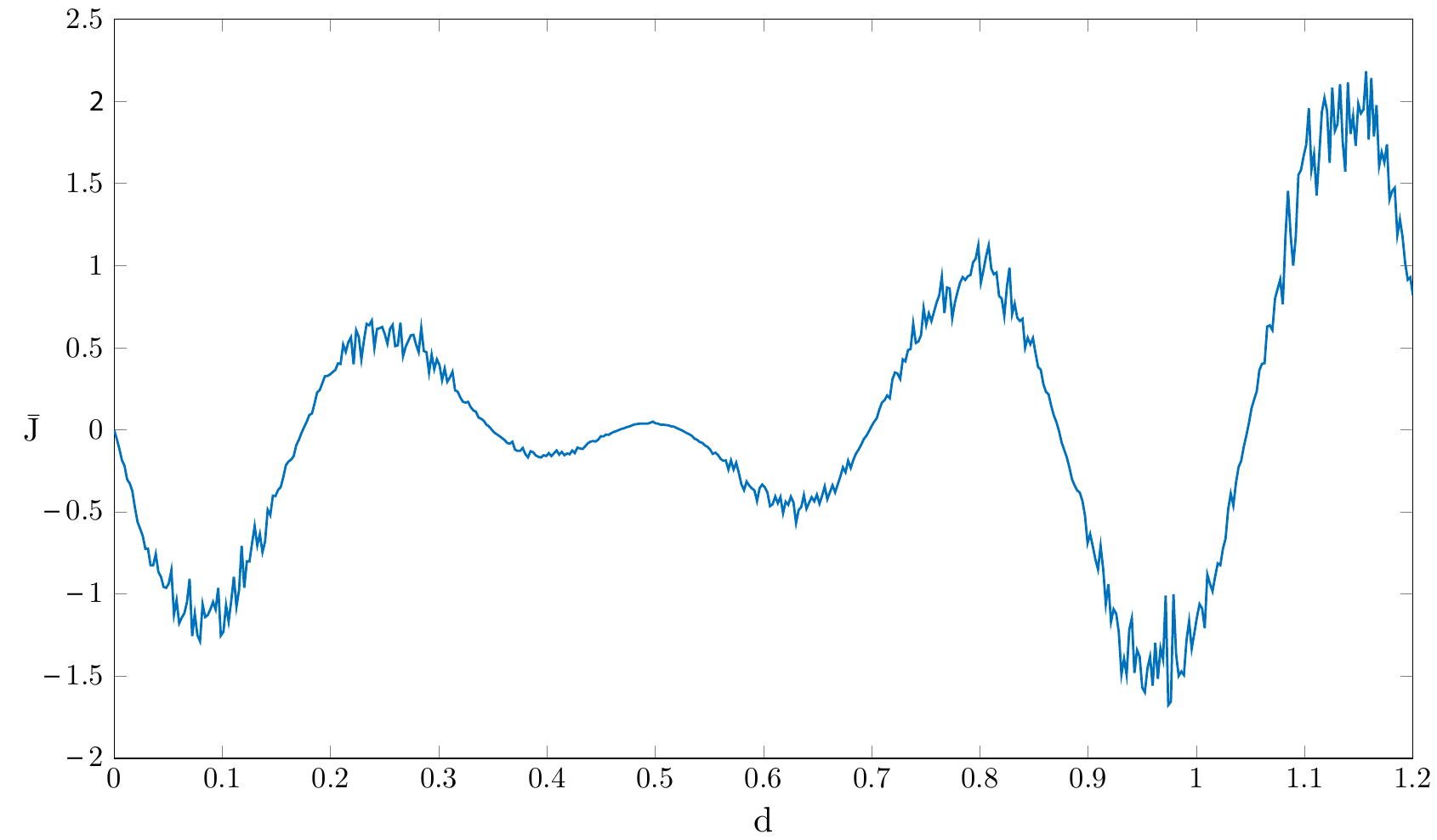}
    \caption{Approximation of $J$ given by MCI}
    \label{fig:noisy_function}
\end{figure}

In order the better visualize the optimization procedure, consider the following multimodal function: 

\begin{equation}
\phi(d, X) = -(1.4-3  d)\sin(18  d) X,
\label{eq:problem_1}
\end{equation}

\noindent where $X$ is a Normal random variable with mean $1$ and standard deviation $0.1$ (\ie $X \sim  \mathcal{N} (1, 0.1)$). Here, we want to find $d^*$ that minimizes the expected value of this function in the design domain $d \in S =  [0, 1.2]$. That is, our problems resumes to:

\begin{equation}
\underset{d \in S}{\min J(d)} = \E\left[\phi(d,X)\right] =  \int_{\Omega} \phi(d,x) \ \pdf{X}(x) \ d x,
\label{eq:problem_1a}
\end{equation}

\noindent where $\E$ is the expected value operator, $\pdf{X}$ is the Normal probability density function of $X$ and $\Omega = (-\infty,+\infty)$ its support. Figure~\ref{fig:noisy_function} illustrates the approximation of $J$ given by MCI over the design domain, \ie it plots $\bar{J}$, evaluated using \autoref{eq:approxmean}, over the design domain.

\begin{figure}[!htbp]
    \centering
    \subfloat[$\overline{\sigma}^2_{target}$ = 1.00]{\includegraphics[width=7.8cm, height=5.6cm]
    {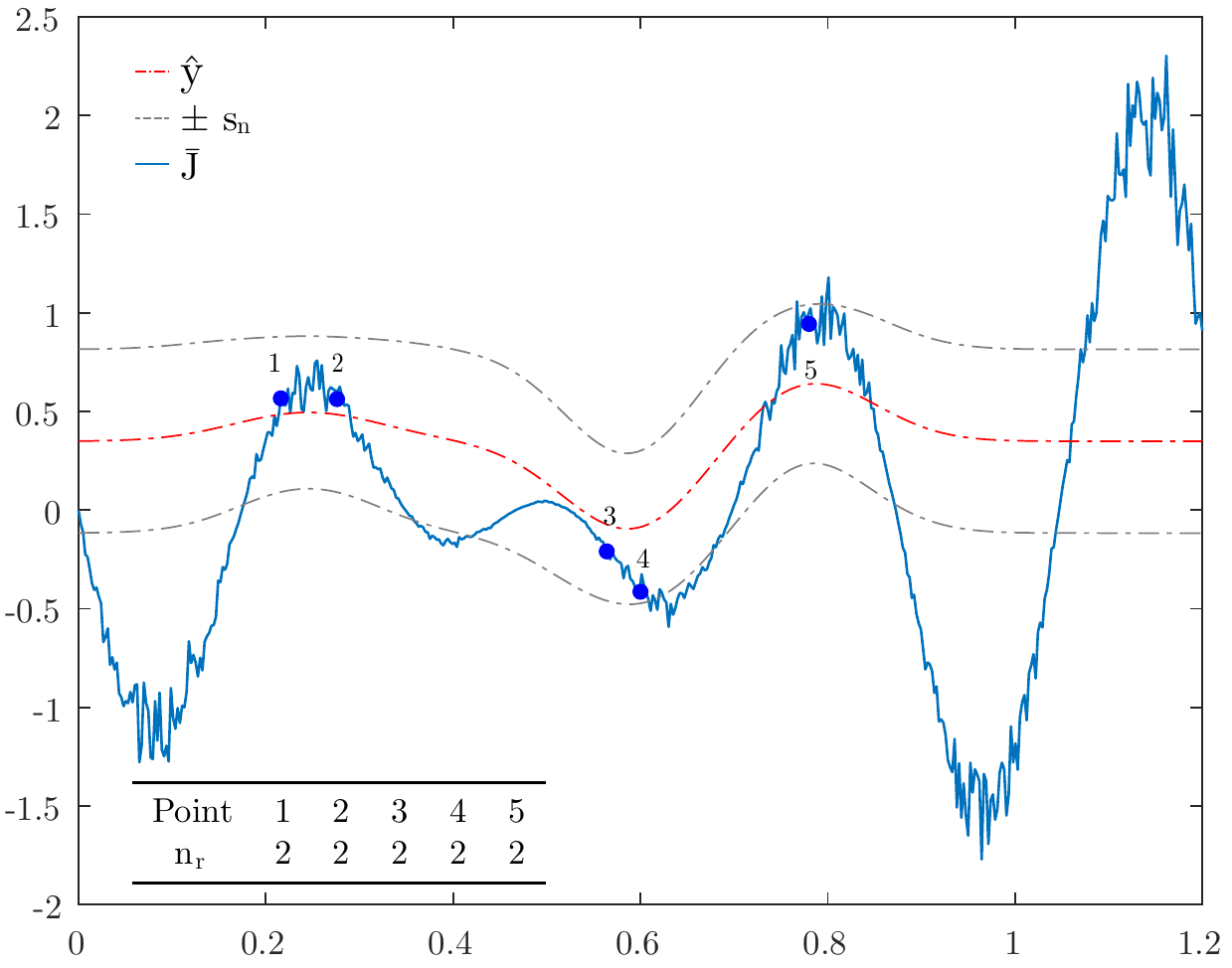} \label{fig:diff_a}}
    \hfill
    \subfloat[$\overline{\sigma}^2_{target}$ = 0.001]{\includegraphics[width=7.8cm, height=5.6cm]
    {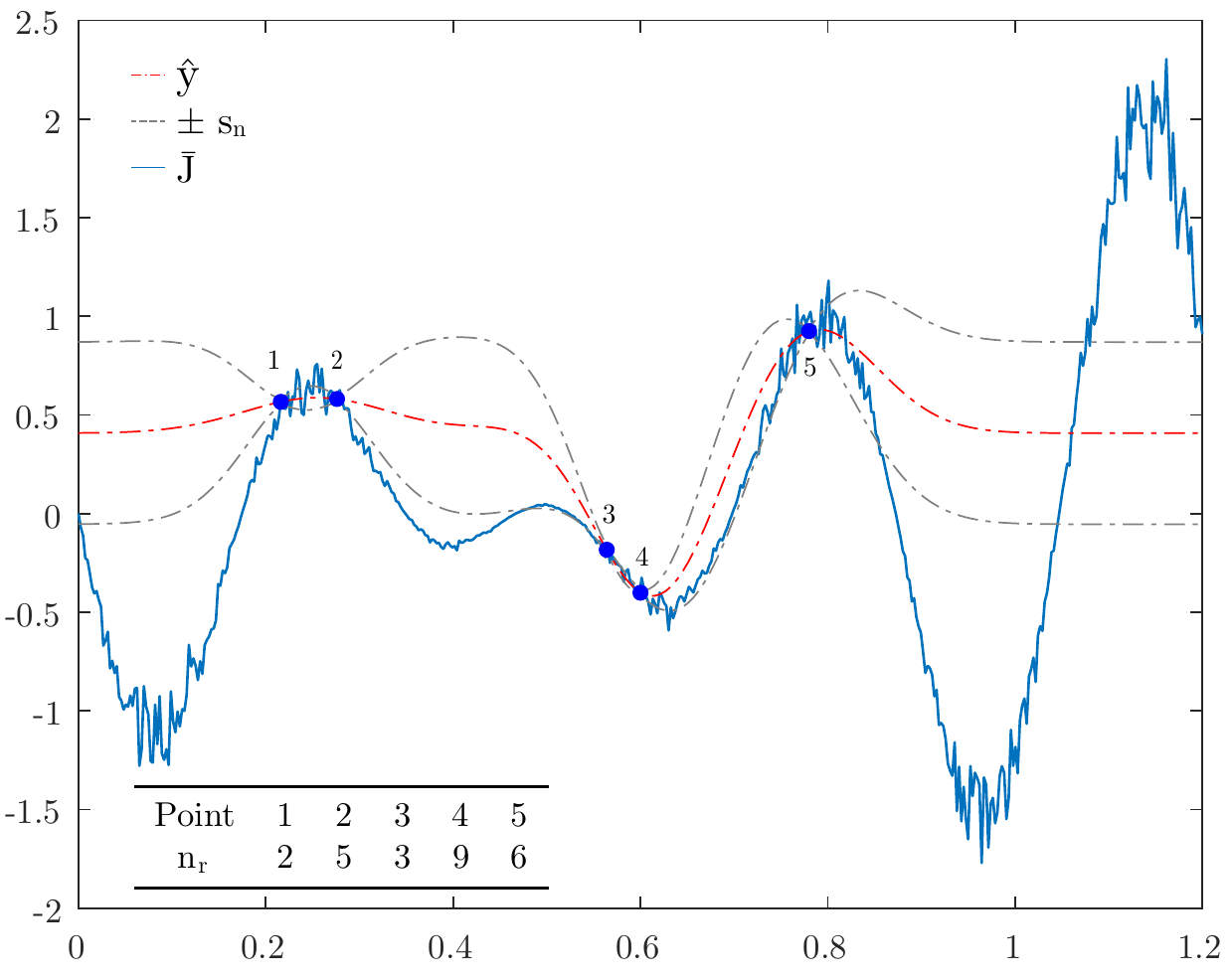} \label{fig:diff_b}}
    \caption{Different model based on the estimated error}
    \label{fig:diff_target_var}
\end{figure}

As described in \autoref{sec:ego}, the first step of an EGO algorithm is the construction of the initial sample plan. Here, we employ the initial Kriging sample size as $n_s = 5$, which may be visualized as the blue dots in Figs. \ref{fig:diff_a} or \ref{fig:diff_b}. Each observation is replicated $n_r$ times in order to achieve the desired target variance ($\overline{\sigma}^{2}_{target}$). 

As already mentioned, different targets imply different error estimates, \ie  the more replications we draw at a given design vector, the lower is the variability of the MCI estimate. For example, Figure~\ref{fig:diff_a} shows a SK model based on $\overline{\sigma}^{2}_{target} =$ 1.00, while Figure~\ref{fig:diff_b} shows the model resulting from the same initial Kriging sample, but in this case with $\overline{\sigma}^{2}_{target} =$ 0.001. In these figures, the dashed red lines represent the SK model prediction, given by \autoref{predSK}, while the dashed gray lines give the variability of the SK model, which correspond to the interval $[-s_n,+s_n]$ from the prediction, which is evaluated by \autoref{eq:skestimatederror}.

%


In deterministic Kriging, the estimated error is exactly zero at
sampled points. However, when $J$ is estimated using MCI, the error only approaches zero
as the number of replications increases (\ie when $n_r \rightarrow \infty$). For instance, comparing both cases in \autoref{fig:diff_target_var}, it can be seen
that the RMSE interval is closer to the sampled points in
Figure~\ref{fig:diff_b} while it is quite far in
Figure~\ref{fig:diff_a}. Although a larger number of evaluations is needed
in order to achieve a lower target variance, it can be seen that it may
improve the prediction. Having a more accurate description of the function
with a lower estimated error on sampled points makes the infill criterion
approach the deterministic case, which is known to be globally convergent.

\begin{figure}[htbp]
    \centering
    \includegraphics[width=\linewidth]{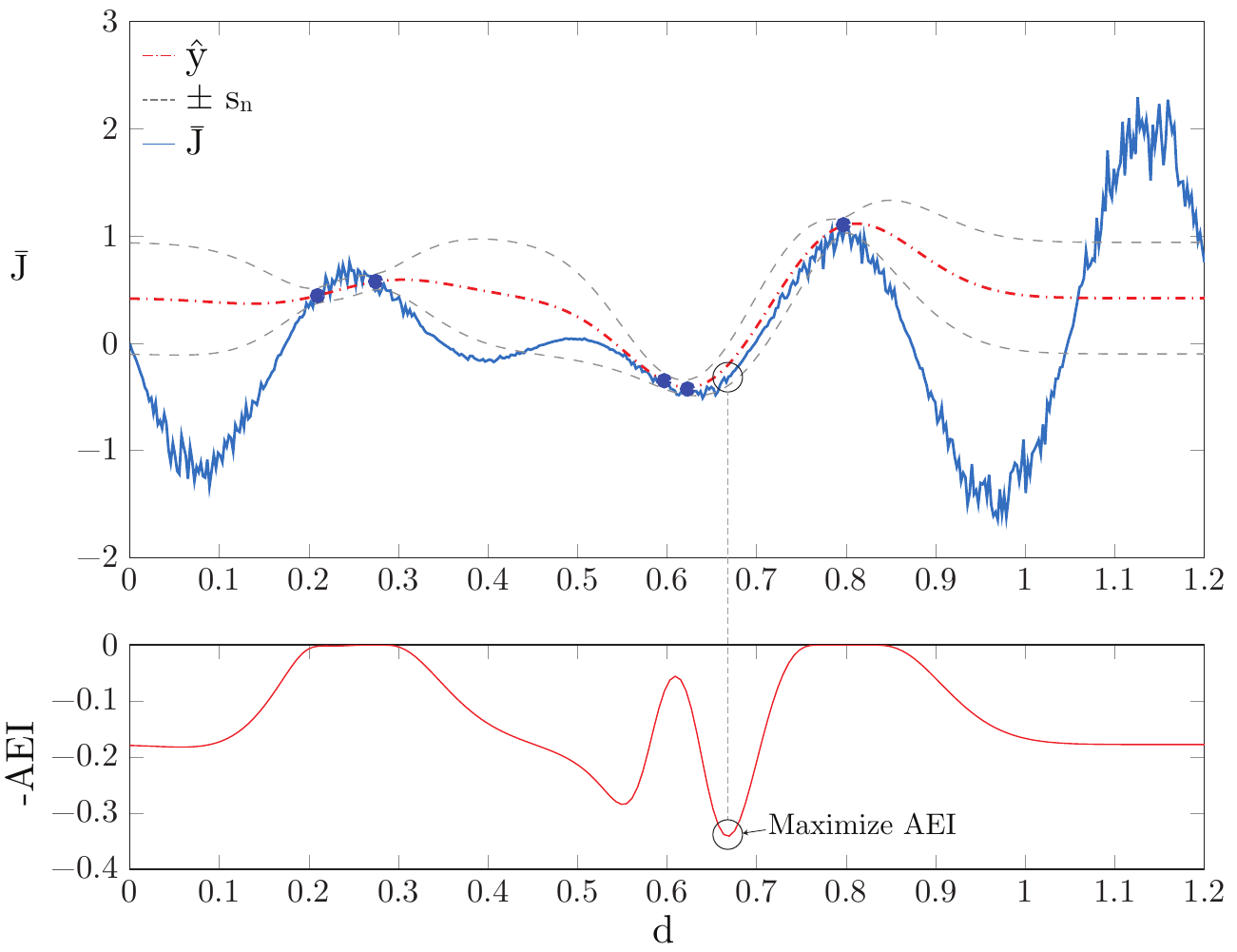}
    \caption{AEI plot over the domain}
    \label{fig:aei_plot}
\end{figure}

After the initial SK model construction, the next step in the EGO
procedure is the refinement of the model by adding infill points. Using the
AEI criterion, the point that minimizes the utility function composed of an
EI term and a penalty term is added to the model (see \autoref{AEIinfill}).

Consider first the case in which we have a relatively low target variance, \ie $\overline{\sigma}^{2}_{target} = 0.01$.
Figure~\ref{fig:aei_plot} illustrates how the AEI infill works. The auxiliary
plot presented below the function plot represents the negative of the AEI
measure over the domain. The circled point indicates the selected infill
point, which is shown to maximize the AEI measure (\ie it minimizes $-AEI$).

\begin{figure}[!htbp]
    \centering
    \subfloat[Infill 5]{\includegraphics[width=0.8\linewidth]
    {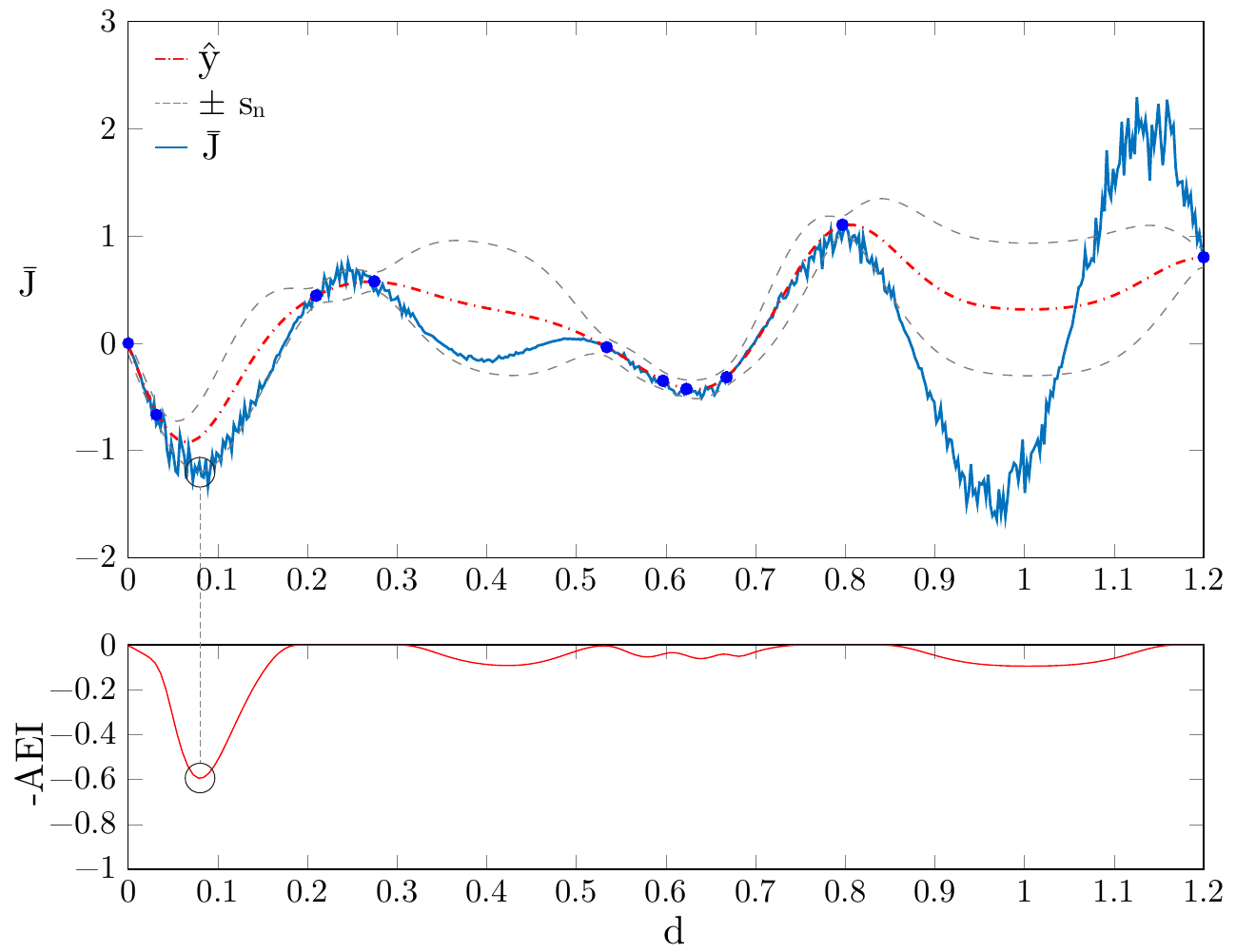}     \label{fig:progress_a}}
    \hfill
    \subfloat[Infill 20]{\includegraphics[width=0.8\linewidth]
    {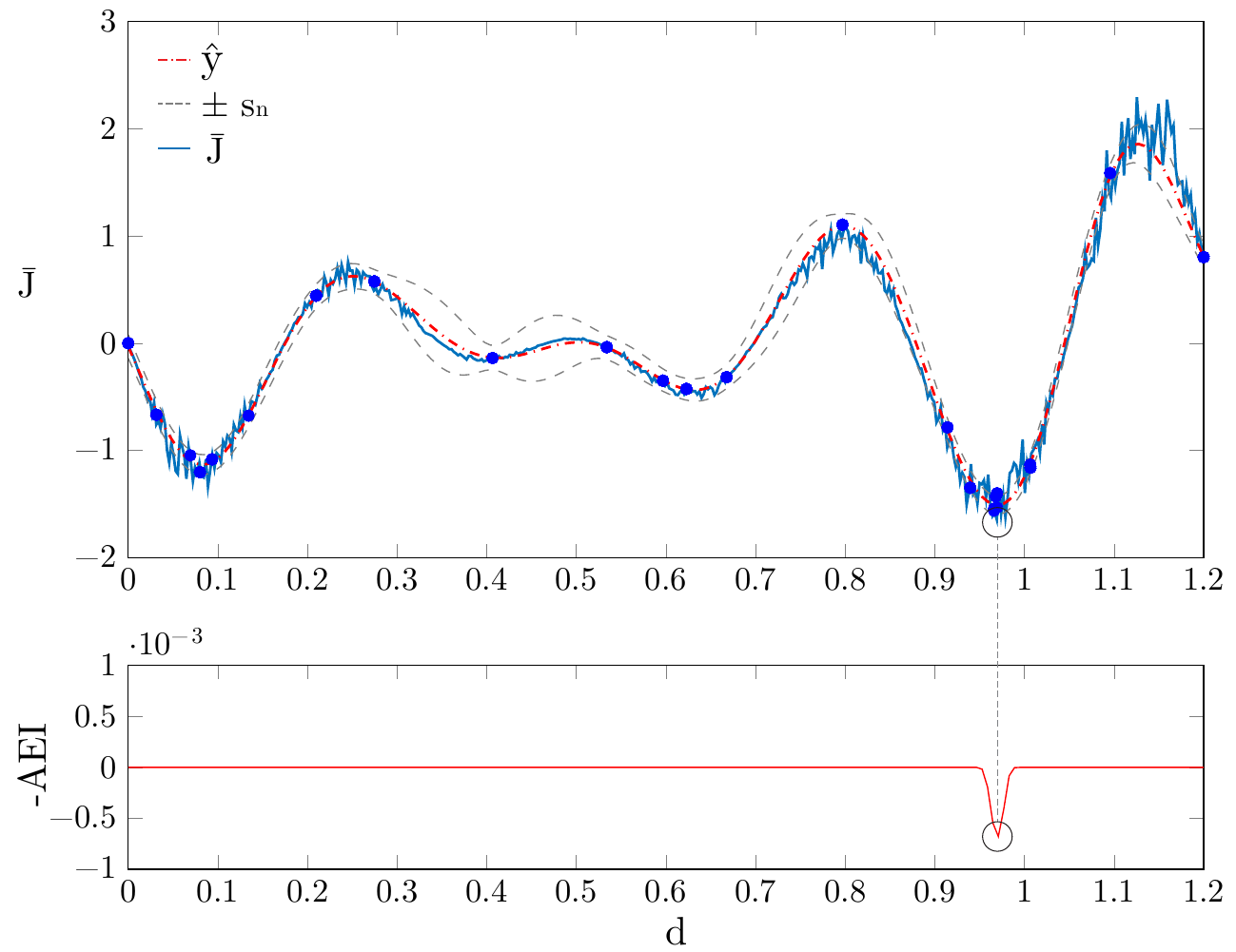}
    \label{fig:progress_b}}
    \caption{Refinement of the SK model seeking the optimal value with target variance $0.01$.}
    \label{fig:progress}
\end{figure}

Figure~\ref{fig:progress} shows the progress of the optimization in two different
moments. In Fig.~\ref{fig:progress_a}, the 5th infill is added to the
model. The search has not found the optimal valley yet. The search proceeds and after
some evaluations on the initial valley the algorithm starts exploring  more
uncertain regions. By the 20th infill, it is possible to see in
Fig.~\ref{fig:progress_b} that  the algorithm already sampled around the
global optimum multiple times. Therefore, AEI value is very low and directed to that single region.

\begin{figure}[htbp]
    \centering
    \includegraphics[width=\linewidth]{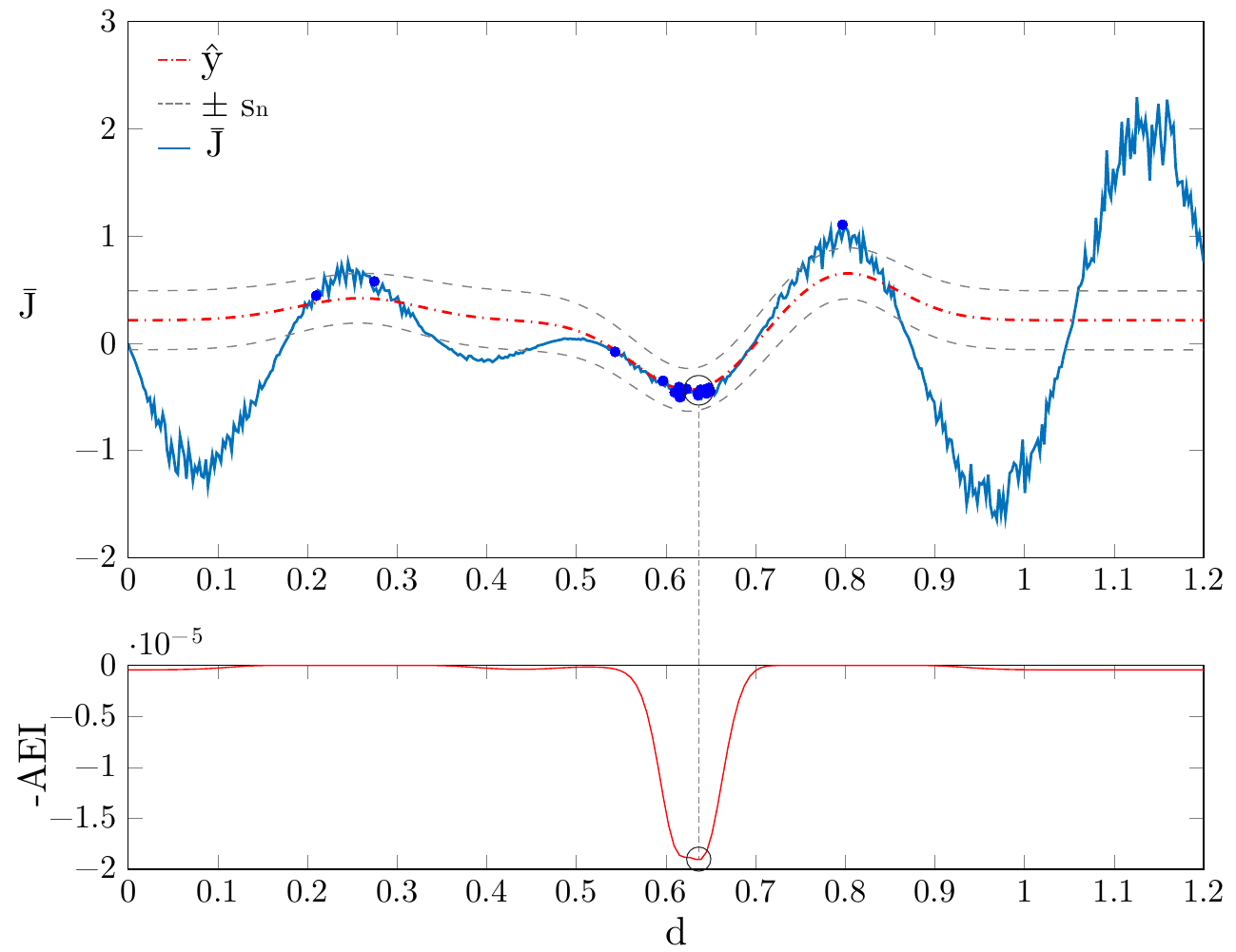}
    \caption{Stalling at infill 20 with target variance 1.00.}
    \label{fig:high_target_stall}
\end{figure}

A completely different situation occurs when we consider the case in which we have a relatively high target variance, \ie $\overline{\sigma}^{2}_{target} = 1.00$. As it may be seen in \autoref{fig:high_target_stall}, which shows the search at the $30^{th}$ infill point, the
optimization process stalls at the first valley found and never explores
other regions. This is caused by the lack of information gained by each
point added to the model. Each infill is inserted considering the unit
target variance. Because the estimated error remains almost the same, the
infill criterion becomes highly local, avoiding exploration of potentially
better regions. The penalization from AEI only takes effect when $s^2_n$ is
low, which is not the case when the target variance is high. 

The main conclusion from this section is that there is
a need for a rational adaptive choice of the target variance. A higher
target reduces computational cost, yet it can become highly local. A lower
target may need higher computational cost, however the refined model is more
easily exploited. A trade-off between these characteristics must be made in
order to:

\begin{itemize}
    \item avoid getting stuck in local minima due to lack of information, 
    \item achieve an efficient optimization process.
\end{itemize}

\section{Proposed approach: adaptive target selection} 
\label{sec:adap_target}

As it was seen in the previous section, the target variance plays an important
role when using the MCI coupled with SK. Setting it either too small or too high may hinder the optimization algorithm. When too small, the target forces a large number of evaluations on a single point, effectively compromising the computational budget of the global optimum search. On the other hand, setting it too high may stall the search due to the lack of reliable information obtained with each infill, as shown in \autoref{fig:high_target_stall}. 


As it can be inferred, the efficiency and the success of the optimization approach
depends on the target selected. However, which target to choose?  The main
contribution of this paper is proposing an adaptive target selection
procedure. This idea arises to confront the existing drawbacks of both low
and high targets, as discussed in \autoref{sec:howtosetvariance}. However,
instead of a simple middle ground constant target variance value, here
we suggest a more efficient approach. 

The proposed approach consists on an adaptive target selection. It aims to
balance the trade-off between accuracy and computational cost. The main idea is
similar to an exploration versus exploitation aspect of any global
optimization procedure. The adaptive approach starts exploring the design domain by evaluating the objective function value of each design point using MCI with a high target variance - so that each evaluation requires only a few samples. That is, the initial target value $\overline{\sigma}^{2}_{target}$ should be initially set to a value that avoids expending too much computational resources at the first iterations of the algorithm, giving room for the adaptive scheme to work. 

Then, it gradually reduces the target variance for the evaluation of additional infill points in regions of the design domain where points were already sampled. It is expected that such a reduction make the algorithm exploit the promising regions of the design domain, not stalling the search as the example of \autoref{sec:howtosetvariance}.

%

\begin{figure*}[ht]
    \centering
    \includegraphics[width=0.7\linewidth]{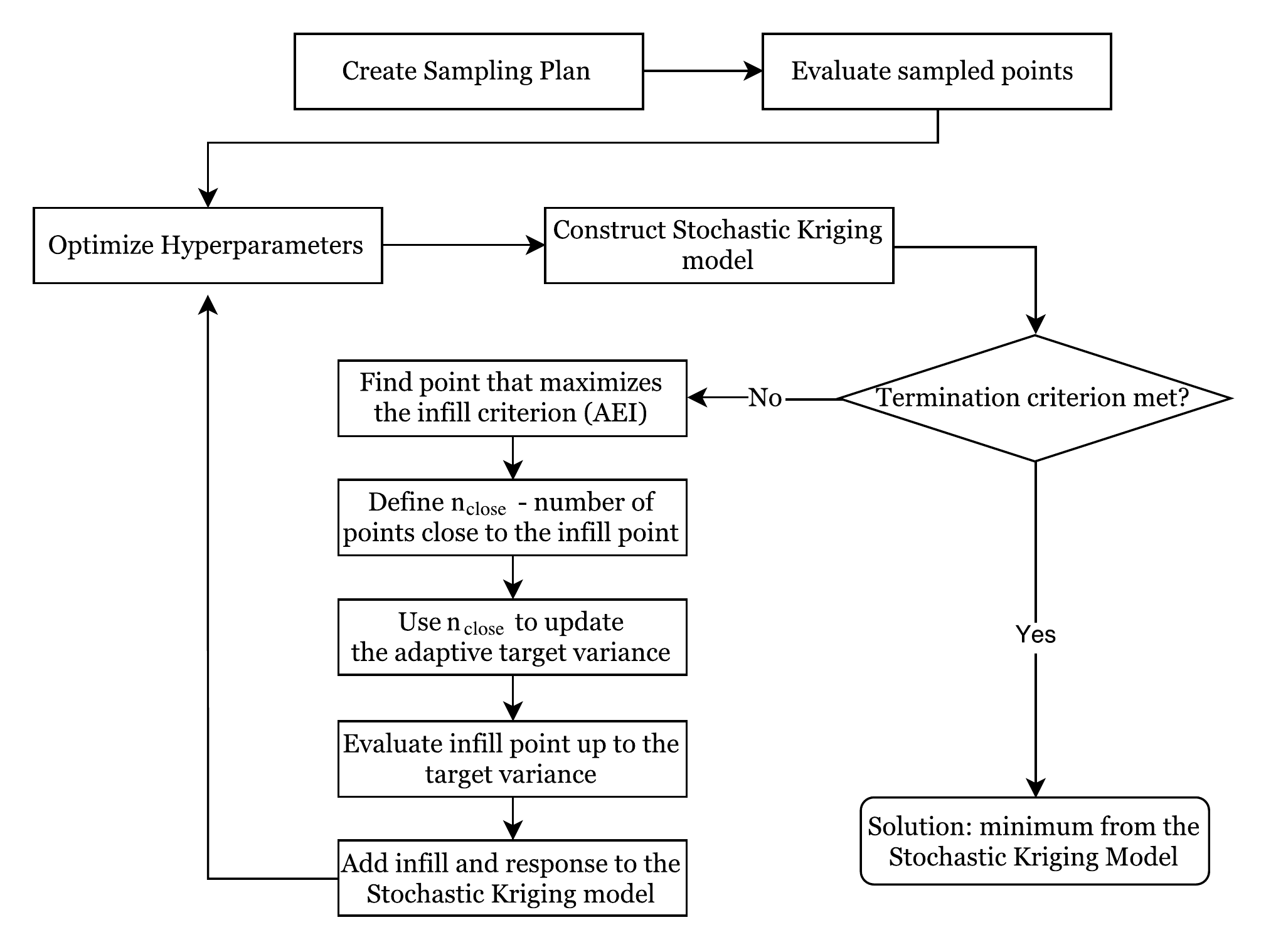}
    \caption{Flowchart of the algorithm}
    \label{fig:fluxogram}
\end{figure*}

A flowchart of the proposed stochastic EGO algorithm, including the proposed
adaptive target selection, is shown in \autoref{fig:fluxogram}. In the next paragraphs, each of its steps is detailed.

The first step is the creation of the Kriging sampling plan. Here, we employed the
Latin Hypercube method presented in \citet{forrester2008engineering}. It is
important to highlight here that all design points of the sampling plan are
evaluated using $n_r = 1$, ignoring the default target variance. The reason
for this is to prevent expending computational resources on the initial
model and in turn, better employ them when promising regions of
exploitation are identified.

After the construction of the SK metamodel for the initial sampling plan, the infill stage begins. The AEI method, presented in \autoref{sec:infill_criteria_for_noisy_evaluations}, is employed for this purpose. Here, an initial target variance $\overline{\sigma}^2_{target}$ is set and the first infill point is added to the model being simulated up to this corresponding target
variance. 

From the second infill point on, the adaptive target selection scheme starts to take place. We propose the use of an exponential decay equation parametrized by problem dimension ($n$) and the number of points already sampled near the new infill point ($n_{close}$).  The latter is defined by the number of points in the model
located at a given distance of the infill point. Suppose \textbf{q} is the point that maximizes the infill criterion. We propose to evaluate $n_{close}$ as

\begin{equation}
n_{close} = \sum_{j=1}^{n_s} I(\mathbf{q},\mathbf{d}^{(j)}),
\label{nclose}
\end{equation}

\noindent where

\begin{equation}
I(\mathbf{q},\mathbf{d}^{(j)}) = \left\lbrace
\begin{array}{ll}
1, & \Vert \mathbf{q} - \mathbf{d}^{(j)} \Vert \leq r_{hc},\\
0, & \Vert \mathbf{q} - \mathbf{d}^{(j)} \Vert > r_{hc},\\
\end{array}
\right.
\end{equation}

\noindent in which $\Vert \ \Vert$ is a given vector norm and $r_{hc}$ is one of the input parameters of the proposed
approach and corresponds to the distance considered around the infill point.

For evaluation of \autoref{nclose} we take, without loss of generality, the maximum vector norm. For an arbitrary vector $\mathbf{u}$, it is given by

\begin{equation}
\Vert \mathbf{u} \Vert = \max_{i = 1,2,..,n} |u_i|,
\end{equation}

\noindent and thus we consider a hypercube around the infill point selected with half-sides $r_{hc}$. Intuitively, this parameter controls the aggressiveness of the search. Increasing it makes the hypercube larger, allowing more sampled points to be treated as close ones.


Then, when the infill is located within an unsampled region, its target variance is set as the initial target variance. On the other hand, when the infill is located in a region with
existing sampled points, a lower target variance ($\overline{\sigma}^{2}_{adapt}$) is employed for the approximation of its objective function value. This is done to
allocate more computational effort on regions that need to be exploited.
Thus, it indicates the purpose of the infill. Isolated infill points focus
on exploring the landscape, where higher MCI accuracy is not needed. When they
start to group up, the focus changes to landscape exploitation. In this
situation, the target MCI variance is set to a lower value, increasing the model
accuracy. By doing so, it also avoids the clustering of multiple
inaccurate points that causes the stalling observed in \autoref{fig:high_target_stall}.

Since different targets are being employed, $\lambda$ becomes
dependent on the sampling characteristics. It is now adaptively updated according to

\begin{equation}
    \lambda(\de) =  
\begin{cases}
\overline{\sigma}^{2}_{target} & \text{ if } \de \in \boldsymbol{\Gamma} \text{ and } n_{close} = 0 \\
\overline{\sigma}^{2}_{adapt}  & \text{ if } \de \in \boldsymbol{\Gamma} \text{ and } n_{close} > 0\\
\overline{\sigma}^{2}_{target} & \text{ if } \de \notin \boldsymbol{\Gamma} \\
\end{cases}, 
\end{equation}

\noindent where $\overline{\sigma}^{2}_{target}$ is the initial target variance and $\overline{\sigma}^{2}_{adapt}$ is the adaptive target variance. The value of $\overline{\sigma}^{2}_{adapt}$ is then set in order to ensure that points in regions already sampled have a low target variance, in order to exploit a promising region accurately.

In this paper, the expression proposed to calculate the adaptive target value for each iteration of the sEGO
algorithm is

\begin{equation}
    \overline{\sigma}^{2}_{adapt} = \overline{\sigma}^{2}_{target}\exp\left(-g\left(n,n_{close}\right)\right),
    \label{eq:target_adapt}
\end{equation}

\noindent where $g$ is a function that depends on the problem dimension ($n$) and $n_{close}$. The reasoning behind the construction of \autoref{eq:target_adapt} is twofold: 

\begin{itemize}
    \item \textit{it displays an exponential decay of the target value}: the choice of an exponential decay seems to be the most intuitive
    considering how the number of function evaluations increases and the error
    decreases. If a linear model were used instead, the target would decrease
    very slightly with few closer points. Yet, it would drop abruptly as
    $n_{close}$ increased, reaching the lower bound for the target. This would cause, initially, an
    unnecessary number of infill points added to the model without a reasonable
    gain of information. Further in the optimization, the target would drop
    abruptly resulting in a large number of evaluations. With the proposed
    approach, the target starts with high values and is progressively lowered
    when points begin to cluster around an optimum valley.
    \item \textit{the decay rate is proportional to the problem dimension $n$}: with
    low dimensional problems, the design space is relatively small so that it becomes
    easier for the infill points to cluster. Thus, the targeting decay cannot
    be too aggressive at risk of expending too much computational resources. At
    higher dimensions, $n_{close}$ does not increase so fast. Thus, it allows
    for a more significant target decay. Hence, it indicates the necessity of $g$ to depend on $n$.
\end{itemize}

\begin{figure}[ht]
    \centering
    \begin{tikzpicture}
    \begin{semilogyaxis}[
        xmin=0,   xmax=16,
        xlabel=$n_{close}$,
        ylabel={$\overline{\sigma}^2$}
    ]
    \addplot+[mark=none, domain=1:15,thick] {1/exp(1/2 + 1/2 * 2 + 9/19 * x  -1/100*2*x)};
    \addplot+[mark=none, domain=1:15,thick] {1/exp(1/2 + 1/2 * 6 + 9/19 * x  -1/100*6*x)};
    \addplot+[mark=none, domain=1:15,thick] {1/exp(1/2 + 1/2 * 10 + 9/19 * x  -1/100*10*x)};
    \addplot[color=gray, dashed, mark=none] coordinates {(0,1e-6)(15,1e-6)};
    \node at (axis cs:1,1e-6) [anchor=south west] {$\overline{\sigma}^{2 \text{ min}}_{adapt} $};
    \addplot [gray, only marks,mark=*] coordinates { (0,2.9e-3) };
    \node [color=gray,pin={[pin distance=2mm]-65:$\overline{\sigma}^2_{target}$}] at (0,2.9e-3) {};
    \legend{$n=2$,$n=6$,$n=10$}
    \end{semilogyaxis}
\end{tikzpicture}
    \caption{Target decay varying with problem dimension}
    \label{fig:decay}
\end{figure}
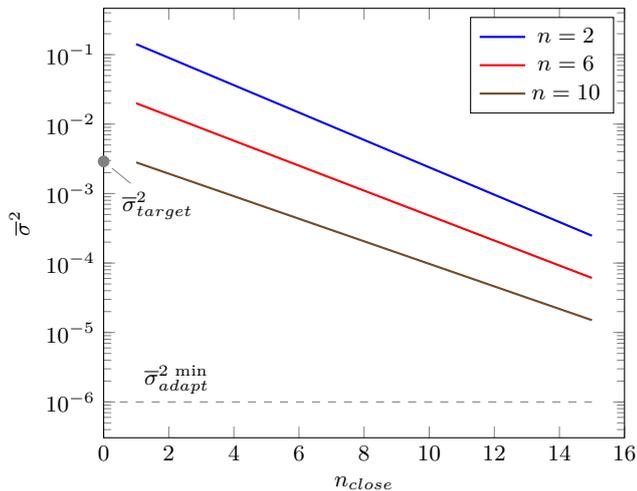

In all the examples of this paper, we employ:

\begin{equation}
g\left(n,n_{close}\right) = a_1 + a_2\cdot n + a_3 \cdot n_{close} - a_4 \cdot n_{close}\cdot n,
\label{eq:target_adaptUSED}
\end{equation}

\noindent where $a_i$ are given constants. Figure~\ref{fig:decay} presents the logarithmic scale plot of \autoref{eq:target_adapt}
for different problem dimensions using $a_1=a_2=a_3=1/2$ and $a_4=1/100$.

It is worth to highlight here that it is also important to set a minimum value for the adaptive target to avoid a computationally intractable number of samples. In other words, not to spend the entire computational budget in only a few points. We thus enforce

\begin{equation}
\overline{\sigma}^{2}_{\min} \leq \overline{\sigma}^{2}_{adapt} \leq \overline{\sigma}^{2}_{target},
\end{equation}

\noindent where $\overline{\sigma}^{2}_{\min}$ is a lower bound on the target.

\section{Numerical Examples}
\label{sec:numerical}

In this section, the performance of the proposed approach is evaluated by the solution of stochastic versions of global optimization benchmark test
functions. These functions possess a single global optimum and most of
them are multimodal, making them good candidates to assess the method's exploration and exploitation
capabilities. These stochastic versions are constructed by inserting multiplicative random variables into the functions. The random variables follow known distributions with specified parameters for each problem.

Here, the efficiency of each algorithm is measured by the number of function evaluations (NFE), \ie the number of times the function $\phi$ is evaluated, which also is employed as the stopping criterion in all the examples. Hence, once the maximum NFE is reached, the algorithm is stopped and the current best design is set as the optimum design of the search. Since the benchmark functions are not expensive to evaluate, it becomes pointless to discuss computational effort as measured in processing time. For this reason, this information is not presented.

It is important to point out that the optimization procedure presented depends on random
quantities. Therefore, the results obtained are not deterministic and may change when the algorithm is run several times. For this reason, when dealing with stochastic algorithms, it is appropriate to
present statistical results over a number of algorithm runs \citep{GOMES201859}. Thus, for each problem, the average as well as the 5 and 95 percentiles of the results found over the set of 100 independent runs are presented. These percentiles are represented by the error bars in the figures. These error bars may be seen as a robustness measure of each algorithm, \ie its ability to provide reasonable results independently of the seed of the random number generator.

Regarding Kriging, for the MLE optimization step, PSO with the default parameters from the Matlab
\citep{MATLAB:2015} implementation was employed. For the AEI infill criterion, $\alpha = 1.00$ was used as suggested by \citet{huang2006global}.

In
\autoref{ssec:optimization_benchmarks}, we compare the constant and the proposed adaptive variance target selection while
\autoref{sec:comparative} aims to compare the proposed scheme against a multi-start algorithm.

\subsection{Stochastic EGO: adaptive against constant target}
\label{ssec:optimization_benchmarks}

In this section, we present the numerical examples in increasing order of complexity, \ie we investigate problems with 1 to 10 design variables. For the proposed adaptive target selection, we employ the framework described in \autoref{fig:fluxogram}, while for the constant target value, $\overline{\sigma}^{2}_{target}$ is set constant throughout the search. In the examples presented in the next subsections, we also evaluate the efficiency of these methods for different input random variables levels by varying their standard deviation. It is important to remark that a different initial sampling plan is employed for each independent run of the algorithm. However, the same initial sampling plans are used for both constant and adaptive targeting to keep a fair comparison between them.

For the upcoming examples, the following parameters are kept constant: initial sampling plans comprised of $n_s = 7n$ points, $r_{hc} = 0.1$, $a_1=a_2=a_3=1/2$ and $a_4=1/100$, $\overline{\sigma}^{2}_{\min} = 10^{-6}$, and for the adaptive scheme, the initial variance target value is set as $\overline{\sigma}^{2}_{target} = 1.0$.

It is worth to highlight that all the parameters of the proposed adaptive approach are kept constant throughout this section. By doing this, we aim to evaluate whether the proposed adaptive approach is able to reach reasonable results without having to tune its parameters.

\subsubsection{Multimodal 1D problem}
\label{subs:1d}

Consider the multimodal 1D problem, given by Eqs. (\ref{eq:problem_1}) and (\ref{eq:problem_1a}), presented in Section \ref{sec:howtosetvariance}. We employ the MCI to approximate the integral in \autoref{eq:problem_1a} and the proposed sEGO scheme to search for $d^*$. Here, we solve this problem for two different standard deviation levels, \ie $\sigma_X = 0.2$ \text{and} $0.3$, and three different values of the stopping criterion, \ie $\text{NFE} = 50, 100, 150$.

\begin{figure}[htbp]
    \centering
    \subfloat[$\sigma_X=0.2$]{
    \includegraphics[width=\linewidth]{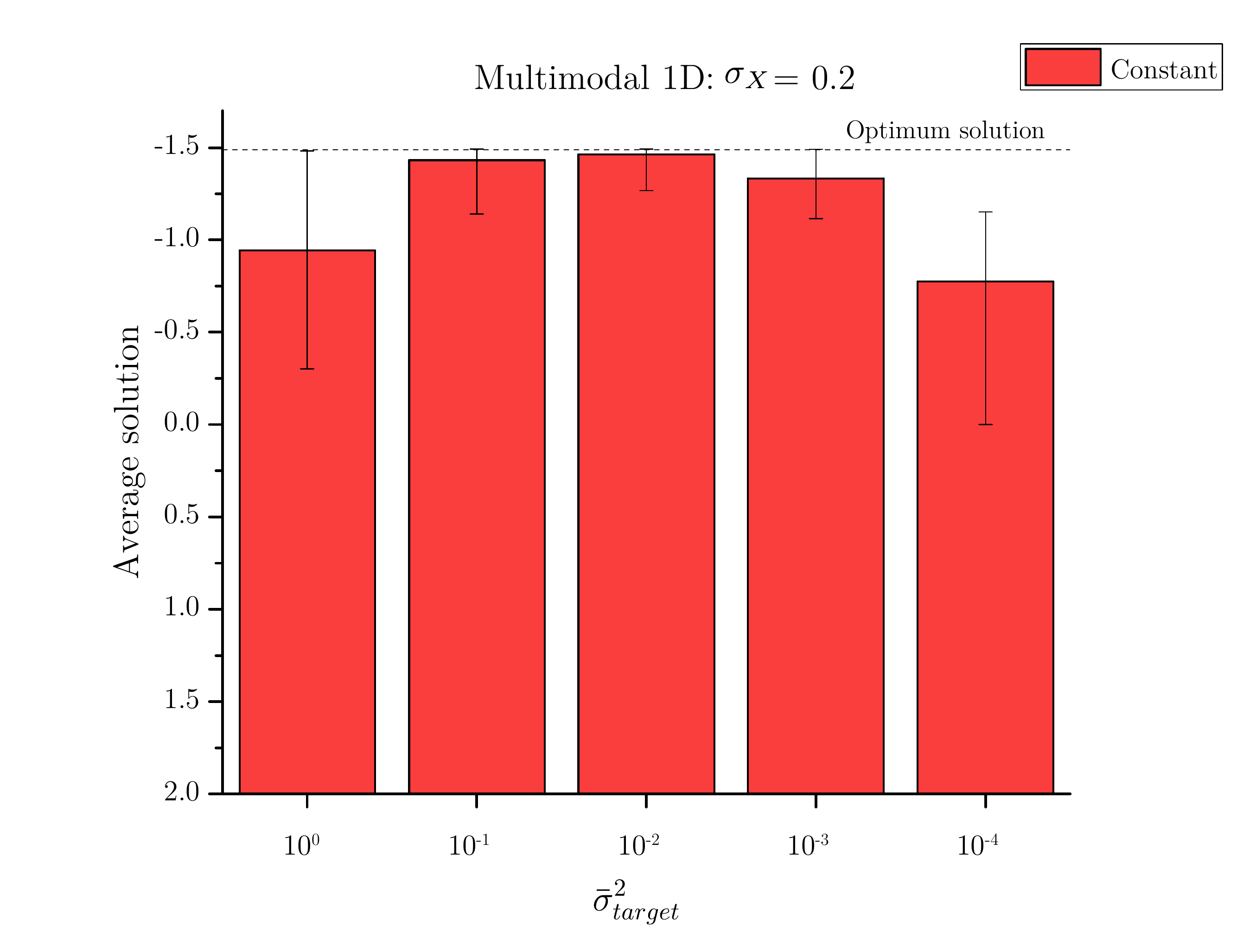}
    }\\
    \subfloat[$\sigma_X=0.3$]{
    \includegraphics[width=\linewidth]{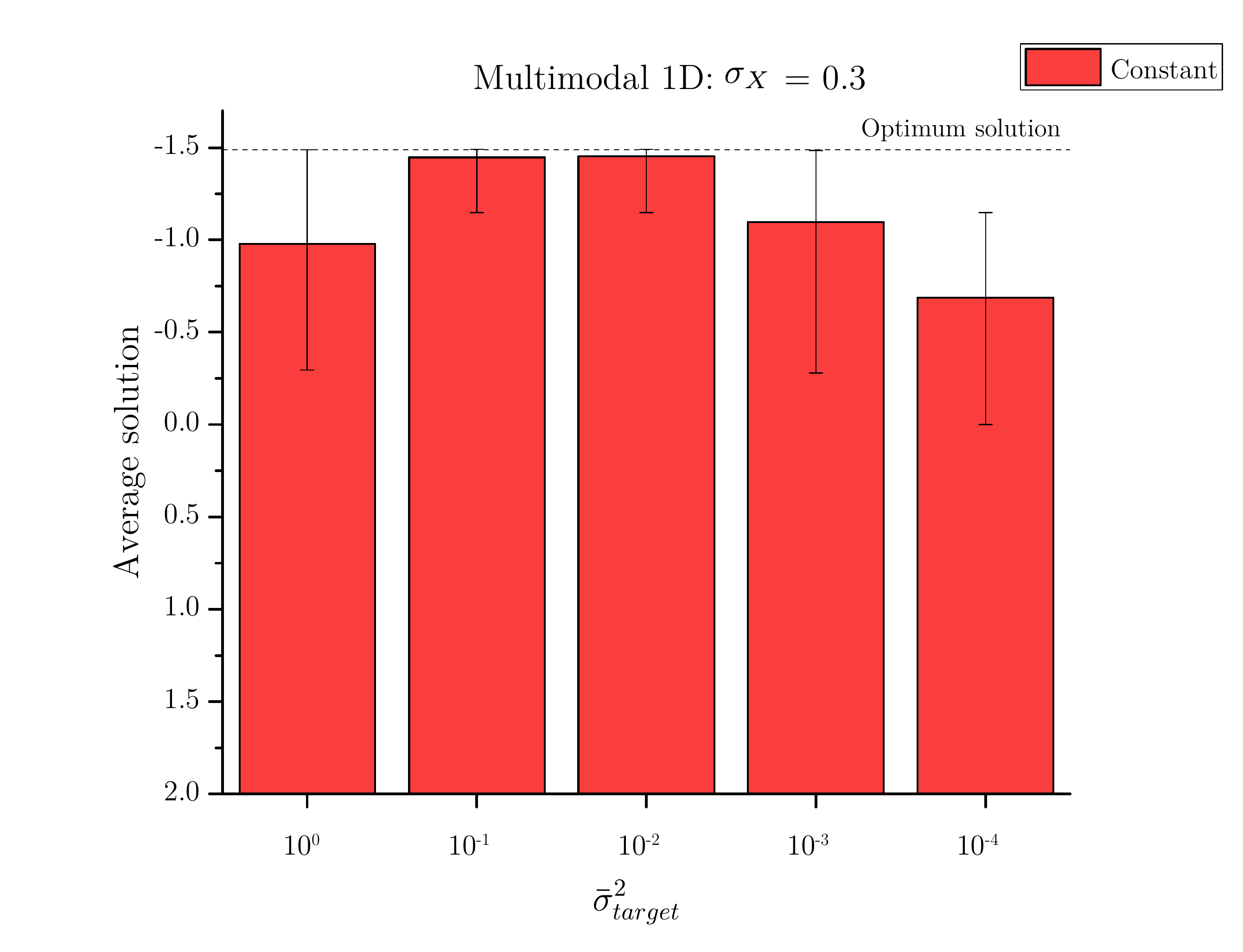}
    }
    \caption{Multimodal 1D problem: results of constant approach for different values of $\overline{\sigma}^{2}_{target}$}
    \label{fig:constant_1d}
\end{figure}

In this example, we first highlight the difficulty of setting a constant target variance. Hence, \autoref{fig:constant_1d} presents the results of the constant approach using different values of $\overline{\sigma}^{2}_{target} = 1.0, 10^{-1}, 10^{-2}, 10^{-3}, 10^{-4}$, and using $\text{NFE} = 150$ as stopping criterion. It is easily noticeable from these results that the performance and robustness of the algorithm is highly dependent on the value of $\overline{\sigma}^{2}_{target}$. For example, for the case $\sigma_X = 0.2$, the best constant variance is $\overline{\sigma}^{2}_{target} = 10^{-2}$ as shown in \autoref{fig:constant_1d}a, while for $\sigma_X = 0.3$, both $\overline{\sigma}^{2}_{target} = 10^{-1}$ and $\overline{\sigma}^{2}_{target} = 10^{-2}$ could be considerd good constant targets, as shown in \autoref{fig:constant_1d}b. It puts in evidence the necessity of an adaptive target variance selection scheme.

\begin{figure}[htbp]
    \centering
    \subfloat[$\sigma_X=0.2$]{
        \includegraphics[width=\linewidth]{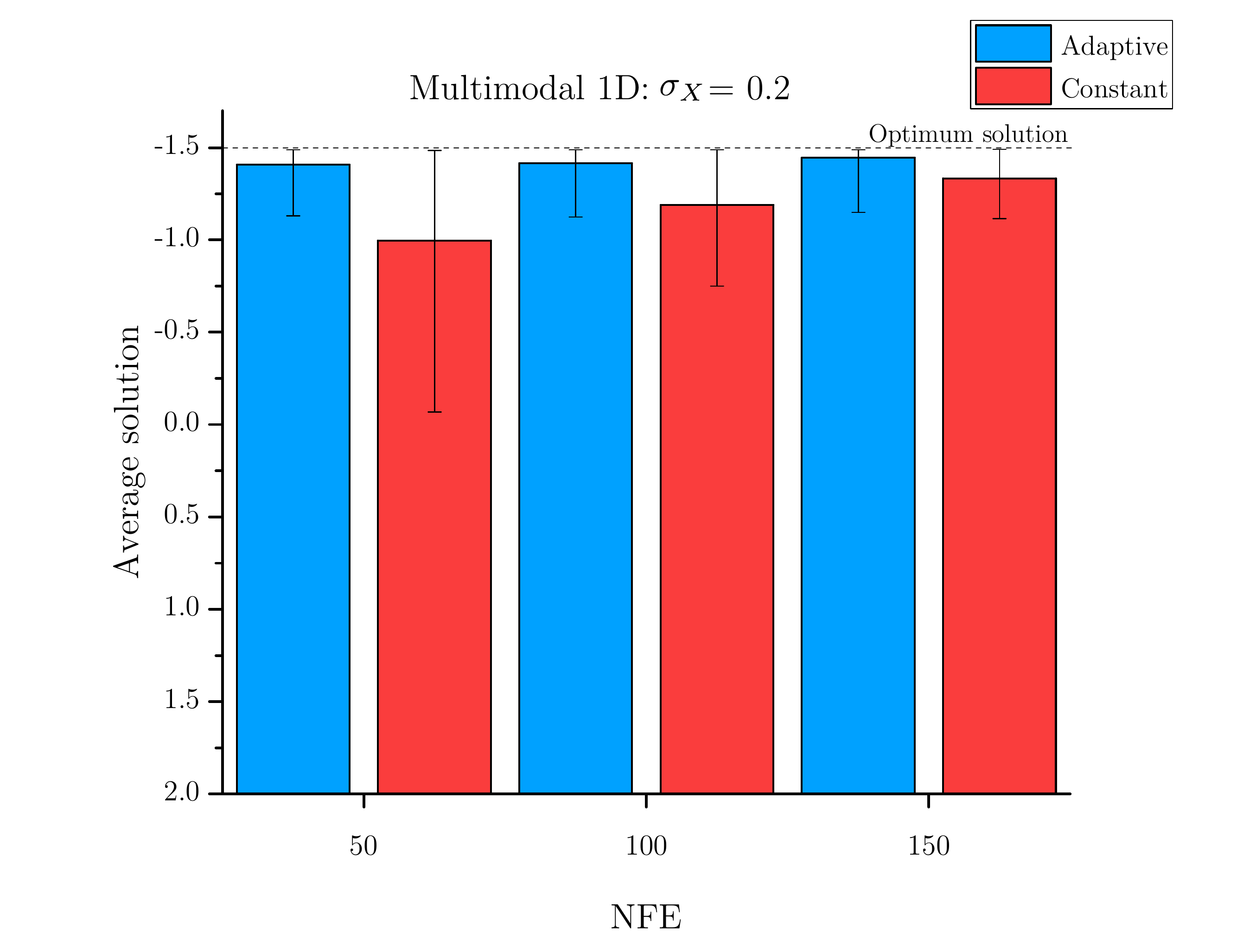}}\\
    \subfloat[$\sigma_X=0.3$]{
        \includegraphics[width=\linewidth]{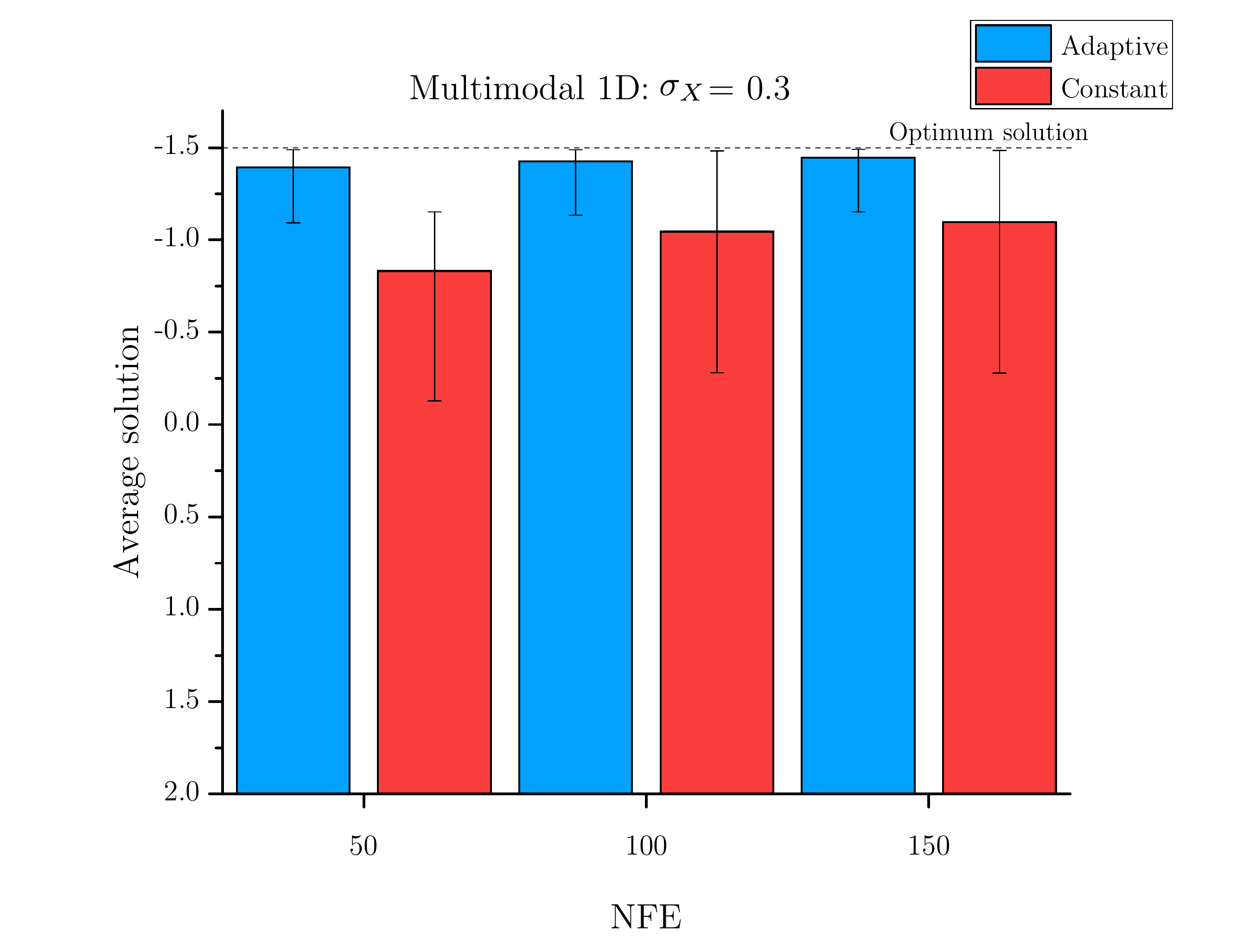} \label{fig:1d_03}}
    \caption{Multimodal 1D problem: results of the adaptive and constant target selection approaches}
    \label{fig:maxeval_1d}
\end{figure}

Figure~\ref{fig:maxeval_1d} presents the results
comparing the constant and adaptive target approaches for different standard deviation values of the input random variable: (a) $\sigma_X = 0.2$ and (b) $\sigma_X = 0.3$. Here, we employ for the constant target scheme,  $\overline{\sigma}^{2}_{target} = 10^{-3}$. In this figure, the height of the bars represent the average solution found over 100 independent runs of the algorithm, while the error bars illustrate the dispersion of the results (5 and 95 percentiles).

For the case in which $\sigma_X = 0.2$, the adaptive approach presented a good performance from the lowest NFE considered, while the constant approach performed better as the computational budget was increased. In the case with higher noise ($\sigma_X = 0.3$), the adaptive approach provided a better value of the mean value of the independent runs as well as lower dispersion. The results presented in this example show that the proposed adaptive target approach successfully minimized \autoref{eq:problem_1a}, and presented clear advantage over the constant approach, especially in the case of higher variability of the input random parameter.

\subsubsection{Multimodal 2D problem}
\label{subs:branin}

In the same manner as the previous subsection, we consider a multimodal
stochastic two dimensional function:

\begin{align}
    \phi(\de, \textbf{X}) = p_1(d_2 -p_2 d_1^2+p_3 d_1-p_4)^2 X_1+
    \nonumber \\
    p_5(1-p_6)\cos(d_1) X_2+p_5+5 d_1,
\end{align}

\noindent with parameters $p_1=1, \ p_2=5.1/(4\pi^2),\  p_3=5/\pi, \ p_4=6, \ p_5=10, \ p_6=1/(8\pi).$
The design domain is $ S =  \{d_1 \in  [-5, 10], \ d_2 \in [0,15]\} $.

$X_1$ and $X_2$ are Normal random variables given by $(X_1, X_2) \sim  \mathcal{N} (1, \sigma_X)$. The weight function $w$ from \autoref{eq:problemabasico} is taken as the probability density function (PDF) $\pdf{\textbf{X}}(\xis)$ of the random vector $\mathbf{X} = \{X_1, X_2\}$. The integral from \autoref{eq:problemabasico} then becomes the expected value of $\phi$ and we have the optimization problem

\begin{equation}
\underset{d \in S}{\min J(\de)} = \E\left[\phi(\de, \textbf{X})\right] = \int_{\Omega} \phi(\de, \textbf{x}) \ \pdf{\textbf{X}}(\textbf{x}) \ d \textbf{x},
\label{eq:IntMin}
\end{equation}

\noindent where $\Omega = \Re^{n_x}$ is the support of $\pdf{\textbf{X}}(\xis)$.

\begin{figure}[!htbp]
\centering
\includegraphics[width=7cm]{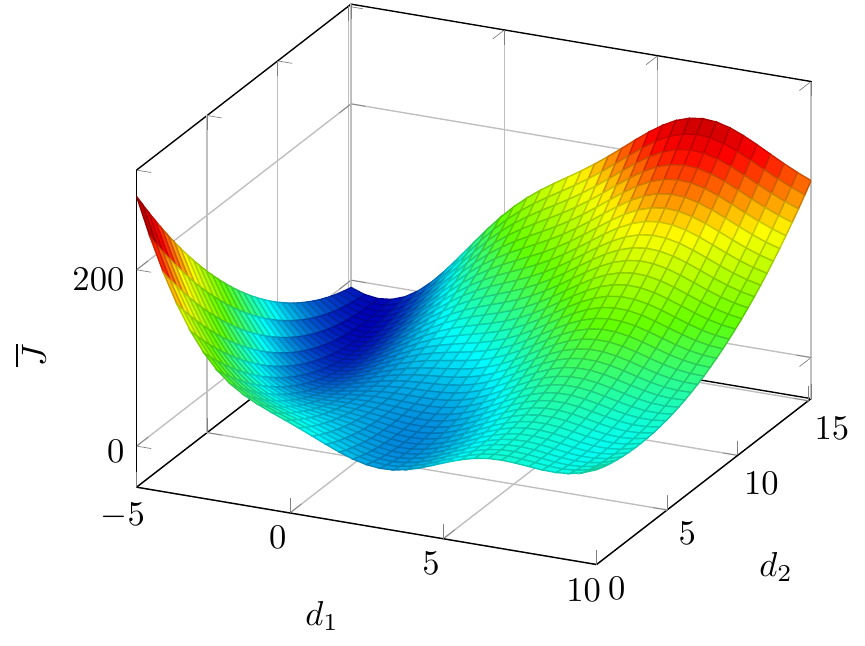}
\caption{Function plot - Branin tilted}
\label{fig:branin}
\end{figure}

In this problem, $\phi(\de,\Xis)$ is a modified Branin function where a term $5d_1$ is added to the
function. This modification forces the existence of a single global optimum
value. The plot of the deterministic version of this function is shown in
\autoref{fig:branin}.


Two cases are considered for different standard deviation values of the input random variable: (a) $\sigma_\mathbf{X} = 0.01$ and (b) $\sigma_\mathbf{X} = 0.05$, and three different values of the stopping
criterion, \ie $\text{NFE} = 100, 300, 1000$. It should be remarked that the behavior of the constant approach is quite similar to the one presented in the previous example, \ie the robustness and performance are highly dependent on the value of $\overline{\sigma}^{2}_{target}$. Hence, we only present from here on the most reasonable results reached by this method. In this example, they are given by $\overline{\sigma}^{2}_{target} = 10^{-2}$.

\begin{figure}[!htbp]
    \centering
    \subfloat[$\sigma_\mathbf{X} = 0.01$]{
    \includegraphics[width=\linewidth]{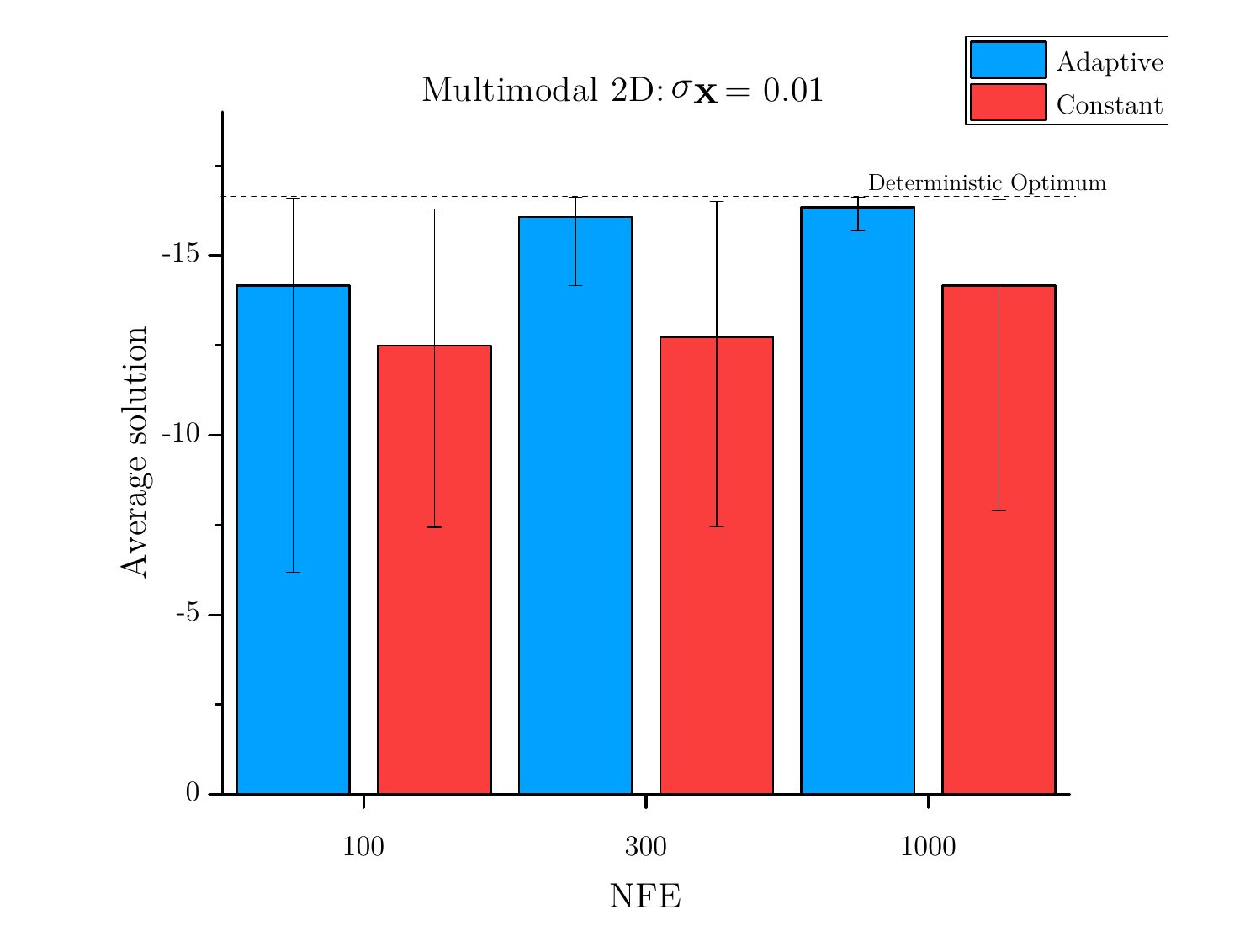}}\\
    \subfloat[$\sigma_\mathbf{X} = 0.05$]{
    \includegraphics[width=\linewidth]{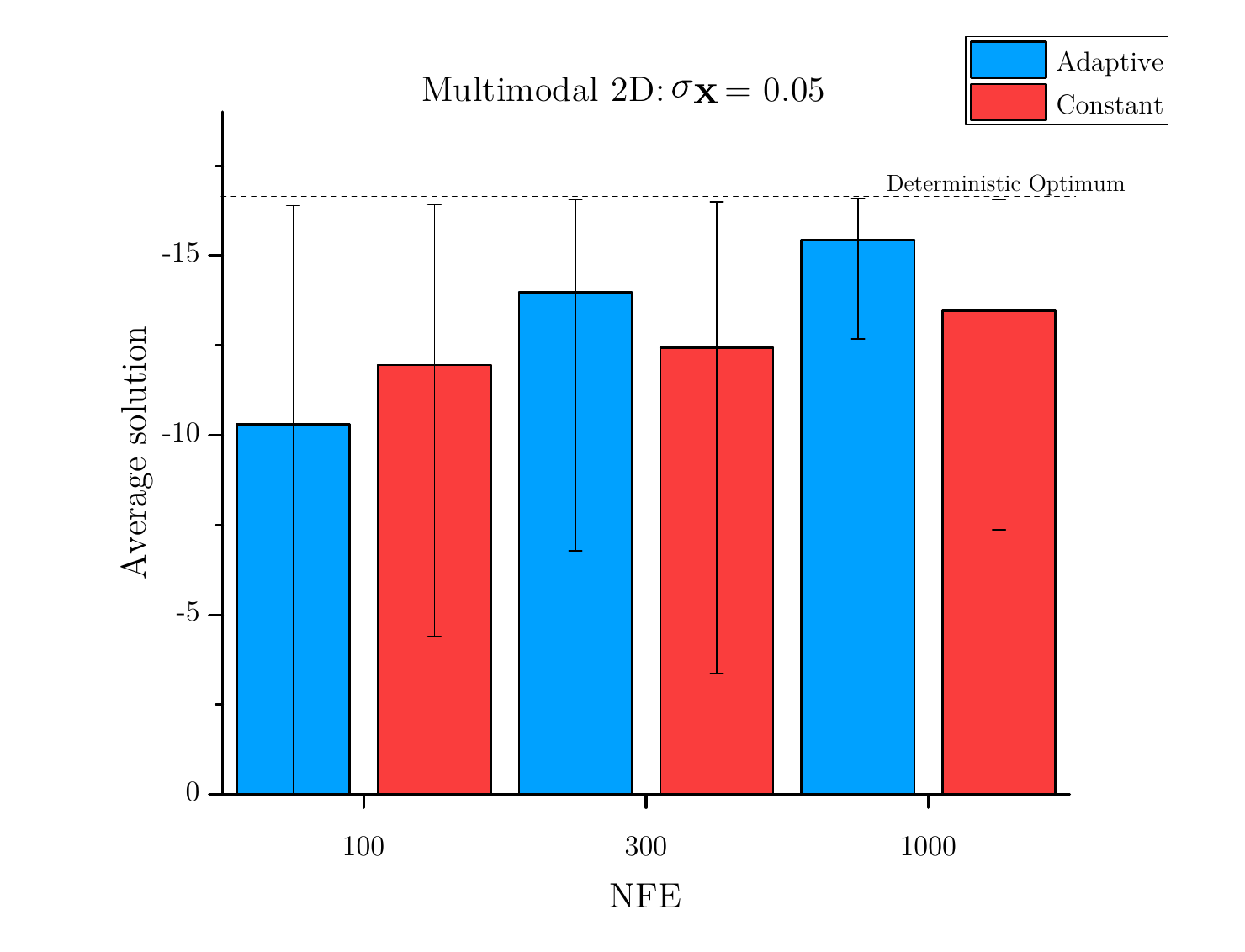}}
    \caption{Multimodal 2D problem: results of the adaptive and constant target selection approaches}
    \label{fig:maxeval_2d}
\end{figure}

Figure~\ref{fig:maxeval_2d} presents the results comparing the constant and
adaptive target approaches. It can be observed that the proposed adaptive approach once more provided better average values and lower dispersion over the independent runs when compared to the constant target approach. Furthermore, in this problem, the adaptive approach provided better results even for the lower noise situation ($\sigma_\mathbf{X} = 0.01$). The only exception was for $\sigma_\mathbf{X} = 0.05$ with NFE=100, where the constant approach presented better results than the adaptive one, as shown in \autoref{fig:maxeval_2d}(b).

\subsubsection{Multimodal 6D problem and high stochastic dimension}

This multimodal benchmark function has six local minima and one global minimum.  It is the well-known Hartman 6D deterministic benchmark function modified by stochastic multiplicative coefficients in the following form:

\begin{equation} {\phi(\de, \textbf{X}) = - \sum_{i=1}^4 {p}_i \exp \left(-\sum_{j=1}^6 A_{ij}
( d_j X_j -P_{ij})^2
\right)},
\end{equation}

\noindent where

\begin{equation}
\mathbf{p} = [1.0, 1.2, 3.0, 3.2]^T,
\end{equation}

\begin{equation}
 A = \left[\begin{array}{rrrrrr}
10   & 3    & 17   & 3.50 & 1.7  & 8 \\
0.05 & 10   & 17   & 0.1  & 8    & 14 \\
3    & 3.5  & 1.7  & 10   & 17   & 8 \\
17   & 8    & 0.05 & 10   & 0.1  & 14
\end{array} \right]
\end{equation}

\noindent and

\begin{equation}
\mathbf{P} = 10^{-4}  \left[\begin{array}{rrrrrr}
1312 & 1696 & 5569 & 124  & 8283 & 5886 \\
2329 & 4135 & 8307 & 3736 & 1004 & 9991 \\
2348 & 1451 & 3522 & 2883 & 3047 & 6650 \\
4047 & 8828 & 8732 & 5743 & 1091 & 381
\end{array} \right].
\end{equation}

The design domain is defined as the unit hypercube $S = \{ d_i
\in  [0, 1], i = 1,2,..,6 \}$. In the following, we analyzed two cases of the Hartman function with different stochastic dimensions: Case 1 with $n_x$ = 6, and Case 2 $n_x$ = 54.

\vspace{.5cm}
\noindent \textit{Case 1: $n_x$ = 6}

The random variables $X_i$ have Normal distribution $X_i \sim  \mathcal{N} (1, \sigma_X)$. The weight function from \autoref{eq:problemabasico} is again taken as the PDF of the random variables and thus the objective function becomes the expected value of $\phi$. The resulting optimization is then as presented in \autoref{eq:IntMin}.

For this problem constant and adaptive approaches are compared for the following standard deviation
values of the input random variable: (a) $\sigma_\mathbf{X} = 0.05$ and (b)
$\sigma_\mathbf{X} = 0.1$. Moreover, three different values of the stopping
criterion were considered: $\text{NFE} = 50, 100, 150$. Finally,
$\overline{\sigma}^{2}_{target} = 10^{-3}$ is employed for the constant
target approach.

\begin{figure}[!htbp]
    \centering
    \subfloat[$\sigma_\mathbf{X} = 0.05$]{
    \includegraphics[width=\linewidth]{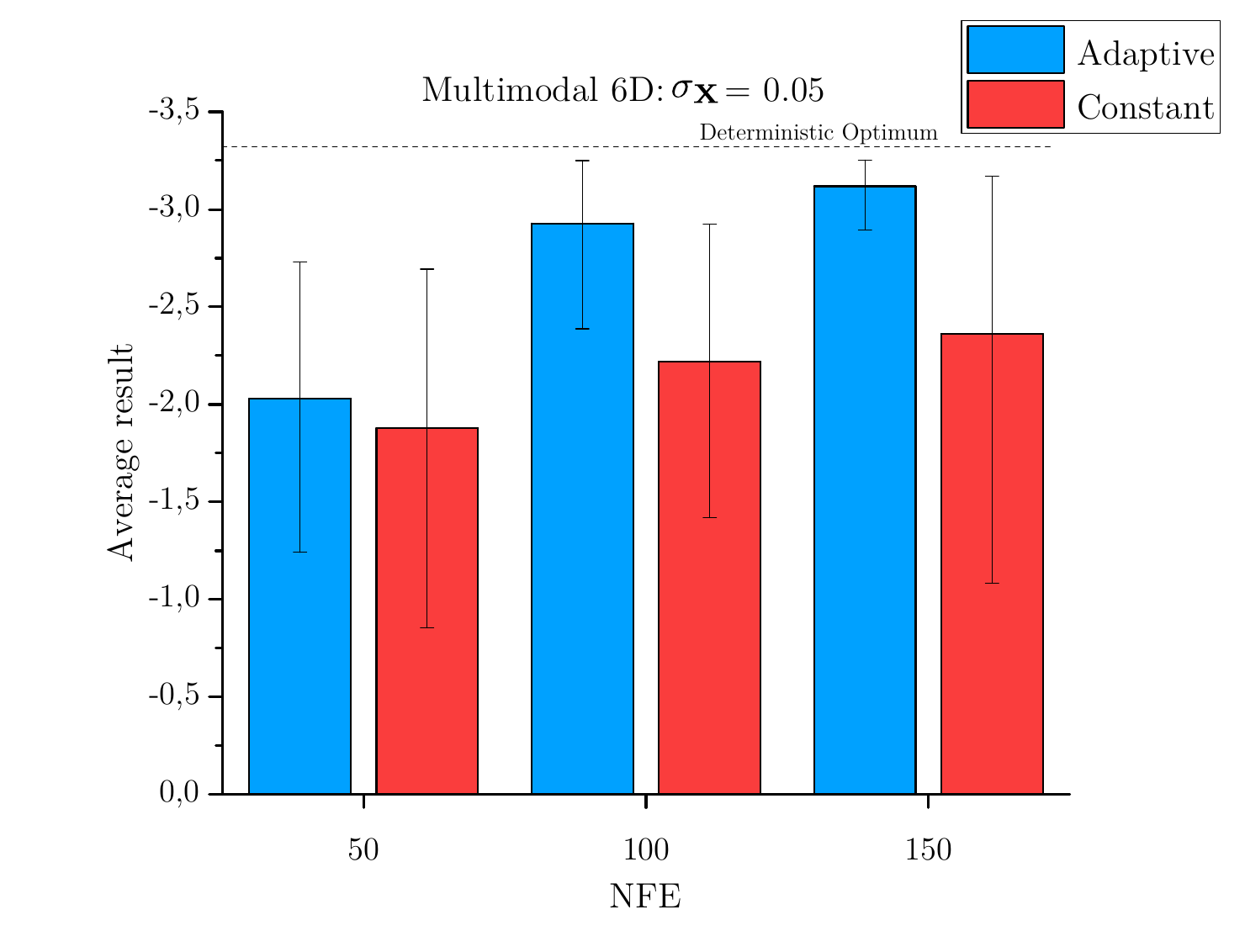}}\\
    \subfloat[$\sigma_\mathbf{X} = 0.1$]{
    \includegraphics[width=\linewidth]{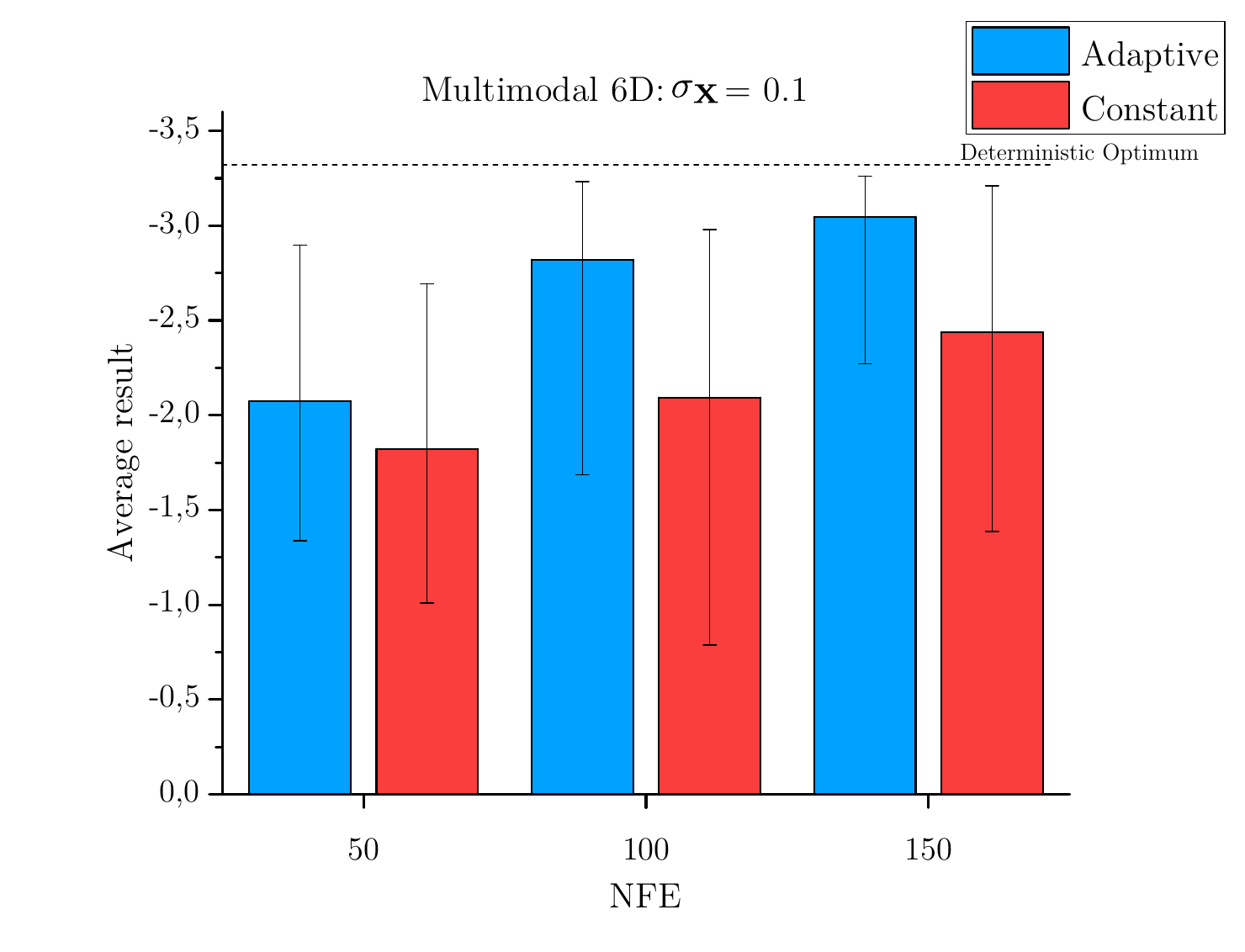}}
    \caption{Multimodal 6D problem: results of the adaptive and constant target selection approaches}
    \label{fig:maxeval_6d}
\end{figure}

Figure~\ref{fig:maxeval_6d} presents the  results for both cases
analyzed. In all the cases, the proposed adaptive approach obtained better results.  As the dispersion of the input random variables increases, it takes more samples to reach a certain variance target, becoming prohibitive to expend the computational budget on accurate exploration points. Consequently, the adaptive target setting enables the optimization to run longer than the constant approach. That is, since the adaptive approach \textit{spends} the available computational budget more rationally, it enables the search to visit more regions on the domain, providing better results. 

Without the early termination proposed in the benchmark, a small constant target would try to reduce the model
uncertainty over the whole design domain instead of only on the promising regions.

\vspace{.5cm}
\noindent \textit{Case 2: high stochastic dimension $n_x$ = 54}

\begin{figure}[!htb]
    \centering
    \includegraphics[width=\linewidth]{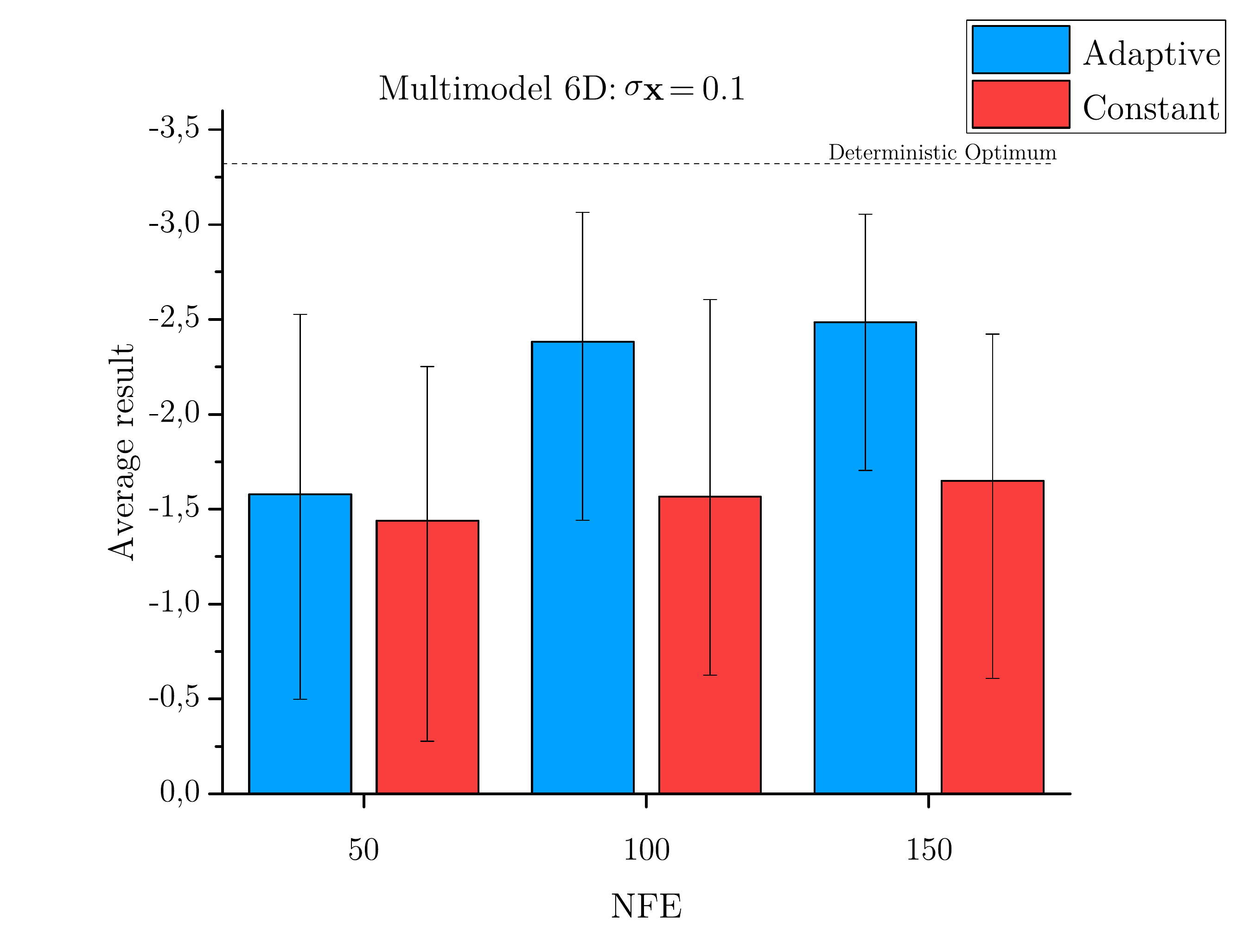}
    \caption{Multimodal 6D problem: higher number of stochastic dimensions}
    \label{fig:6d_highdim}
\end{figure}

This case aims at investigating the performance of the proposed approach in problems with high stochastic dimension ($n_x$). For this purpose, the six random variables of case 1 are kept and we also model all the elements of matrices $\mathbf{A}$ and $\mathbf{P}$ as independent normal random variables. More specifically, $\left(\mathbf{A}\right)_{ij} = A_{ij} \cdot \mathcal{N} (1, \sigma_A)$ and  $\left(\mathbf{P}\right)_{ij} = P_{ij} \cdot \mathcal{N} (1, \sigma_P)$, where $\sigma_A$ and $\sigma_P$ are both 0.01.

Figure~\ref{fig:6d_highdim} presents the results from the analysis of Case 2.
Worse results than those obtained from Case 1 are expected as we are dealing with a higher number of stochastic variables.
Moreover, the same number of evaluations was employed for the stopping criterion, in order to enable direct comparison. Nevertheless, the
adapative procedure presented better results than the constant case. Additionally, variability of the adaptive results decreased with the increase of NFE, showing a convergence pattern. This case illustrates the promising use of the proposed adaptive approach for problems with a high number of stochastic dimensions.

\subsubsection{Multimodal 10D problem}

The last benchmark problem is a stochastic version of the Levy function, which is the following $n$-dimensional multimodal benchmark problem: 

\begin{align}
    \phi(\de, \textbf{X}) &= \sin^2(\pi p_1) + \sum_{i=1}^{n-1}(p_i-1)^2[1+10 \sin^2(\pi p_i+1)]
     \nonumber\\& +(p_n-1)^2[1+\sin^2(2 \pi p_n)], 
\end{align}

\noindent where $p_i = 1 +\frac{d_i X_i}{4}$ for $i=1,2,...,n$. Here we take $n= 10$ and a design domain $S =\{
d_i \in  [-10, 10], i = 1,2,...,10 \}$.  The random variables $X_i$ follow a Normal distribution with $\sigma_\mathbf{X}=0.01$, \ie, $X_i \sim  \mathcal{N} (1, 0.01)$. As we did in the previous examples, the weight function is taken as the PDF of the random vector and the problem is written as the minimization of the expected value of $\phi$, as presented in \autoref{eq:IntMin}. Here, three different values of the stopping criterion were considered: $\text{NFE} = 50, 100, 150$, and we set $\overline{\sigma}^{2}_{target} = 10^{-2}$ for the constant target approach.

\begin{figure}[htbp]
    \centering
    \includegraphics[width=\linewidth]{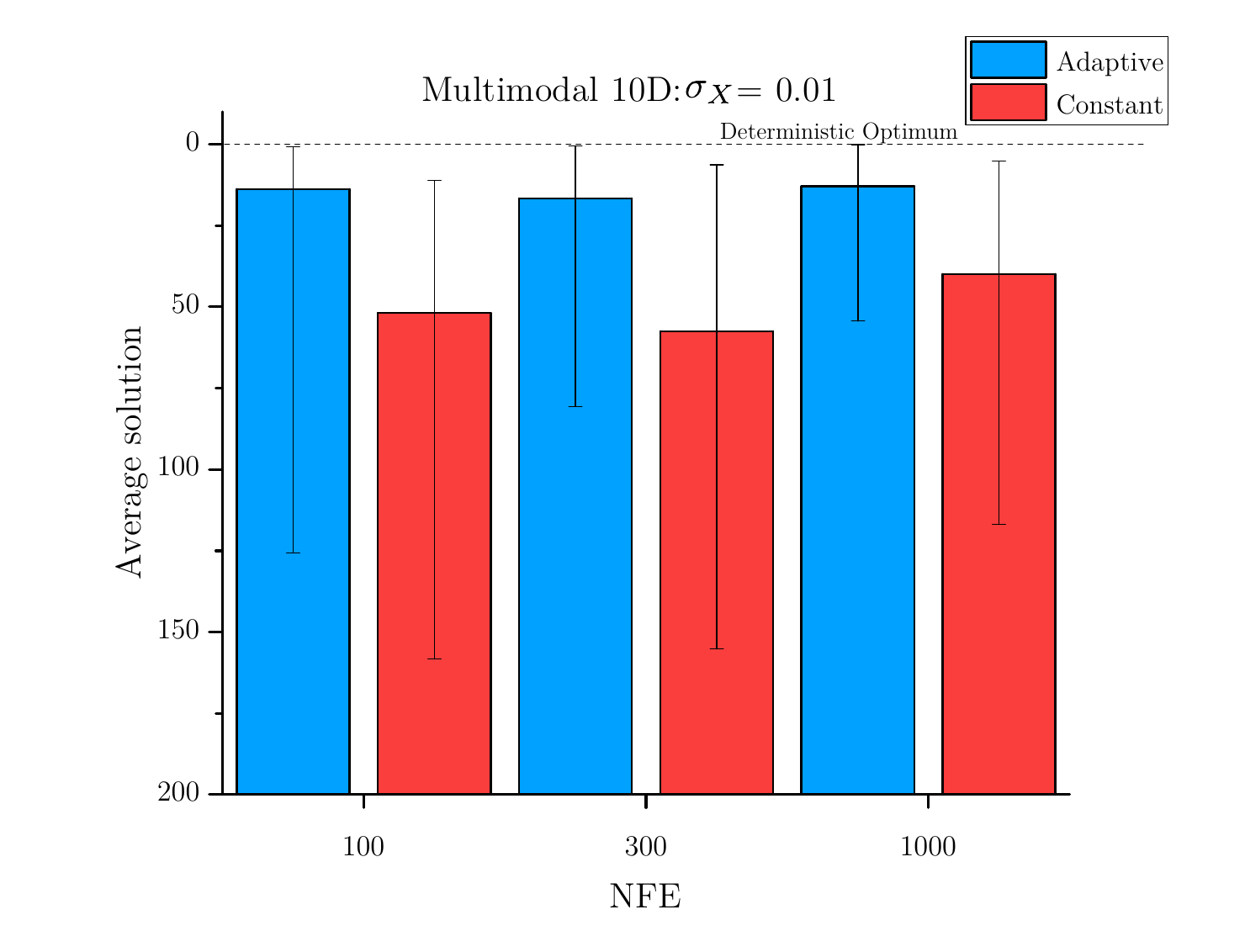}
    \caption{Multimodal 10D problem: results of the adaptive and constant target selection approaches}
    \label{fig:maxeval_10d}
\end{figure}

The results are presented in \autoref{fig:maxeval_10d}.  Similarly to the previous example, the adaptive targeting obtains better
average results than the constant counterpart. Moreover, increasing the
maximum number of function evaluations consistently decreases the
variability of results, represented by the error bars. The method obtain
reasonable results using a relatively small number of function evaluations
even considering a very large 10-dimensional design space.

\subsection{Adaptive approach against multiple start optimization methods}
\label{sec:comparative}

It became clear from the previous examples that the proposed adaptive sEGO algorithm 
can be successfully employed in the optimization of problems of the sort of \autoref{eq:problemabasico}. However, when
investigating the performance of global optimization algorithms, it is
desirable to present a comparison with multi-start based approaches 
\citep{LeRiche2012}. Hence, in this section, we make a comparison with the Globalized Bounded Nelder–Mead (GBNM) algorithm \citep{luersen2004globalized}, which  consists of a probabilistic restart
procedure coupled with a Nelder-Mead local search \citep{nelder1965simplex}. The GBNM was successfully applied to several multimodal problems in engineering, such as laminated composite structures \citep{Luersen2004}, truss optimization \citep{TORII2011}, rotor robust design \citep{Ritto11,Lopez2014}, damage identification \citep{FADELMIGUEL20134241,NHAMAGE201647}, among others. Further details on the algorithm are presented in Appendix~\ref{sec:gbnm}.

\subsubsection{Adaptive target and GBNM with fixed size sample}

\begin{figure}[htbp]
    \centering
    \includegraphics[width=\linewidth]{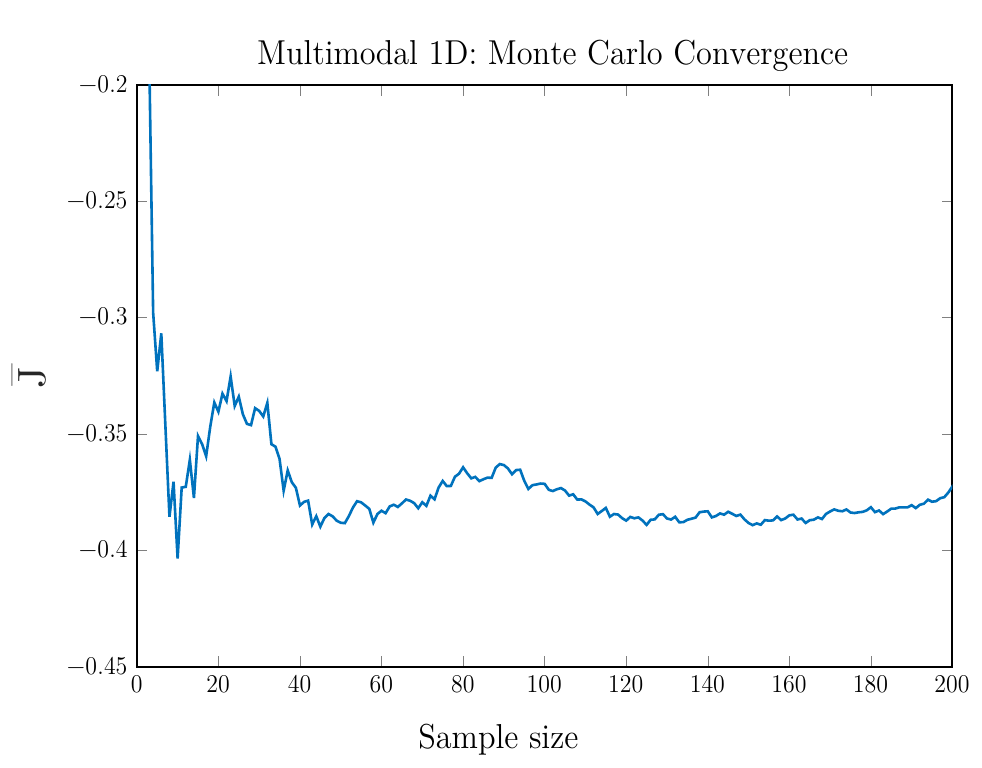}
    \caption{Mean convergence for Multimodal 1D function}
    \label{fig:convergence}
\end{figure}

In order to set the sample size in Eq. (\ref{eq:approxmean}) for the GBNM, we employ the following procedure widely adopted in the literature of robust design \citep{Capiez08a,Capiez08b,Ritto11,Lopez2014,leticia2016,MIGUEL2016703}. First, we construct a convergence plot of $\bar{J}$ with respect to the sample size employed in MCI. Then, we set $n_r$ as the sample number around the value that $\bar{J}$ becomes stable. For example, \autoref{fig:convergence} shows an example of the convergence curve of $\bar{J}$  in a given point of the design domain of the multimodal 1D problem. $\bar{J}$ becomes stable after approximately 100 simulations is employed in MCI. Hence, we set $n_r = 100$, keeping it constant throughout the search, and sample $\xis^{(i)}$ always using  the same seed of the random number generator.

Regarding computational cost, however, note that this is the same multimodal 1D function that we optimized in section \ref{subs:1d}, whose results are in \autoref{fig:maxeval_1d}.
If 100 simulations were employed for each new point required by the GBNM algorithm, the computational cost would be
much higher than the NFE = 150 stopping criterion employed with the proposed approach. Consequently, the proposed approach is much more efficient than any algorithm similar to the GBNM in this problem if this fixed size sample method is employed to obtain $\bar{J}$. Furthermore, the same conclusion may be drawn if we apply this procedure to the other benchmark problems analyzed in Section 5.1. Hence, this analysis is not presented for them here.

\subsubsection{sEGO and GBNM with constant target variance}

\begin{figure}[htbp]
    \centering
    \includegraphics[width=\linewidth]{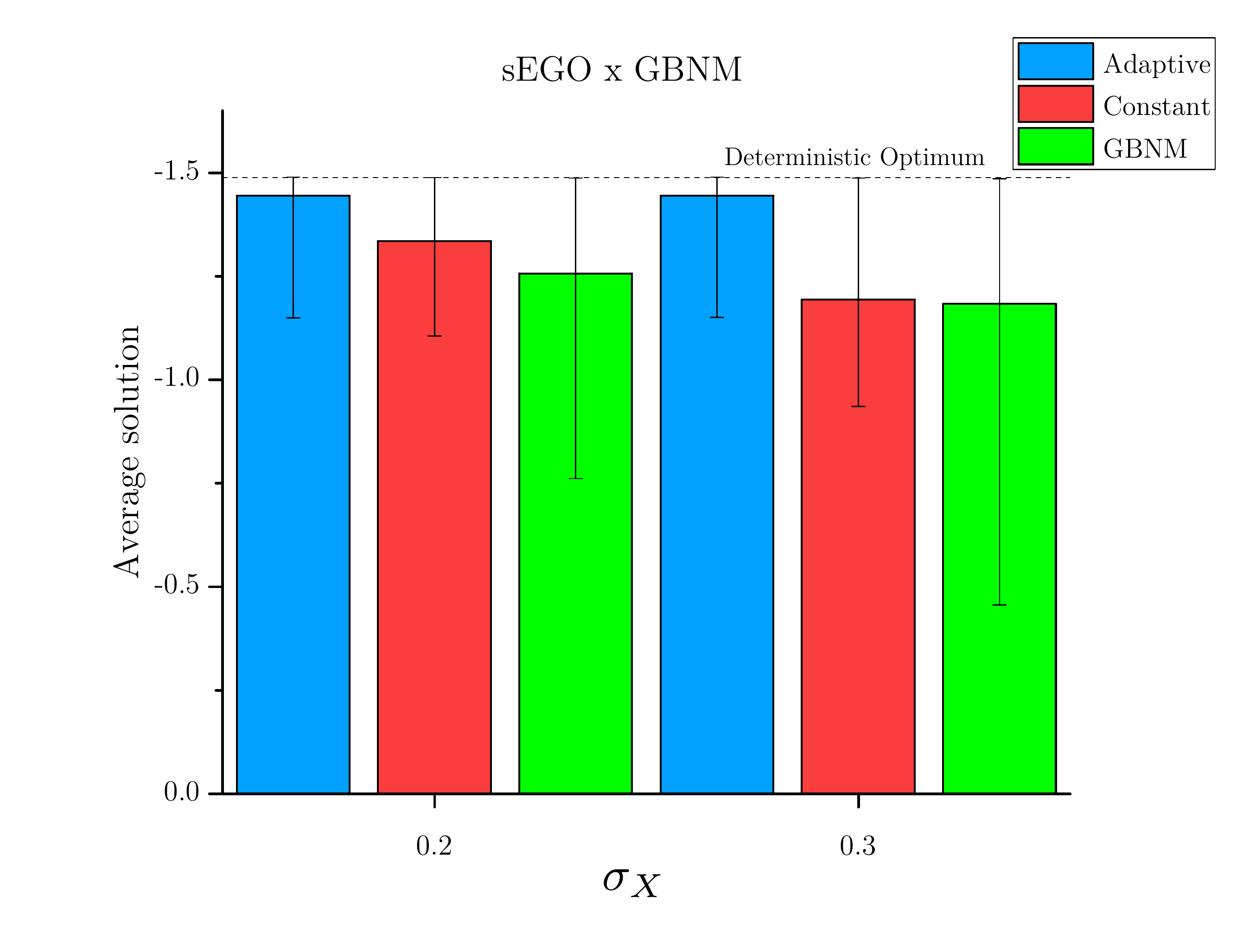}
    \caption{Comparison between sEGO and GBNM for the Multimodal 1D function}
    \label{fig:gbnm_1d}
\end{figure}

Consider again the 1D problem from section \ref{subs:1d}. We now set the same target variance for both algorithms, GBNM and constant approach. That is, we set $n_r$ as the sample size that provides $\overline{\sigma}^{2}_{target} = 10^{-3}$, and the stopping criterion as the maximum number of NFE like the previous section. Two cases are analyzed: (a) $\sigma_X = 0.2$ and (b) $\sigma_X = 0.3$, and \autoref{fig:gbnm_1d} presents the obtained results for the sEGO approaches and GBNM.

In case (a) ($\sigma_X = 0.2$), it can be seen that the sEGO algorithm display a lower dispersion of the results as well as lower objective function value. For the $\sigma_X = 0.3$ case, the maximum number of evaluations differs in each method. For sEGO, we maintain
NFE = 150, while for GBNM, we increase NFE to 1000. This is done because if we set NFE = 150, GBNM does not provide reasonable results. By comparing both approaches in this case (\autoref{fig:gbnm_1d}), it can be seen that GBNM
provides even more spread results, while sEGO obtain results closer to the optimum solution, especially for the proposed adaptive approach.

This behavior persists and is largely amplified when the problem dimension
increases. This leads to an even higher discrepancy  between the sEGO 
approaches (constant and adaptive target) and GBNM. Thus, further examples in higher dimensions are not
shown. The results of this section indicate that the sEGO approaches largely outperform multi-start based methods for problems of the sort of \autoref{eq:problemabasico}.

\subsection{Application: Tuned-Mass Dumper (TMD) system optimization}
\label{sec:tmd_optimzation}

As discussed in \autoref{ModelSituation}, there are several engineering problems where the
proposed integral minimization could be employed. In this section, one application in the field of
structural dynamics is shown. It involves the seismic vibration control of a structure.  The optimization problem resumes to determine the stiffness and damping of a TMD device in order to maximize the structural reliability of a building subject to seismic excitation \citep{chakraborty2011reliability}. The problem is briefly presented in the next paragraphs. For a more detailed description, the reader is referred to \citet{Lopez2015aa}. It should be noted that the problem presented here is strictly academic. Aspects of the real structure are not
taken into account and simplifications were made such as: assuming the structural response linear
elastic, the excitation comes from a stationary process, the random
variables of each floor are uncorrelated, \etc Nevertheless, it remains as a fairly complex problem
comprehending the engineering fields of optimization, control, dynamics and reliability.

\begin{figure}[hbpt]
    \centering
    \includegraphics[width=8cm]{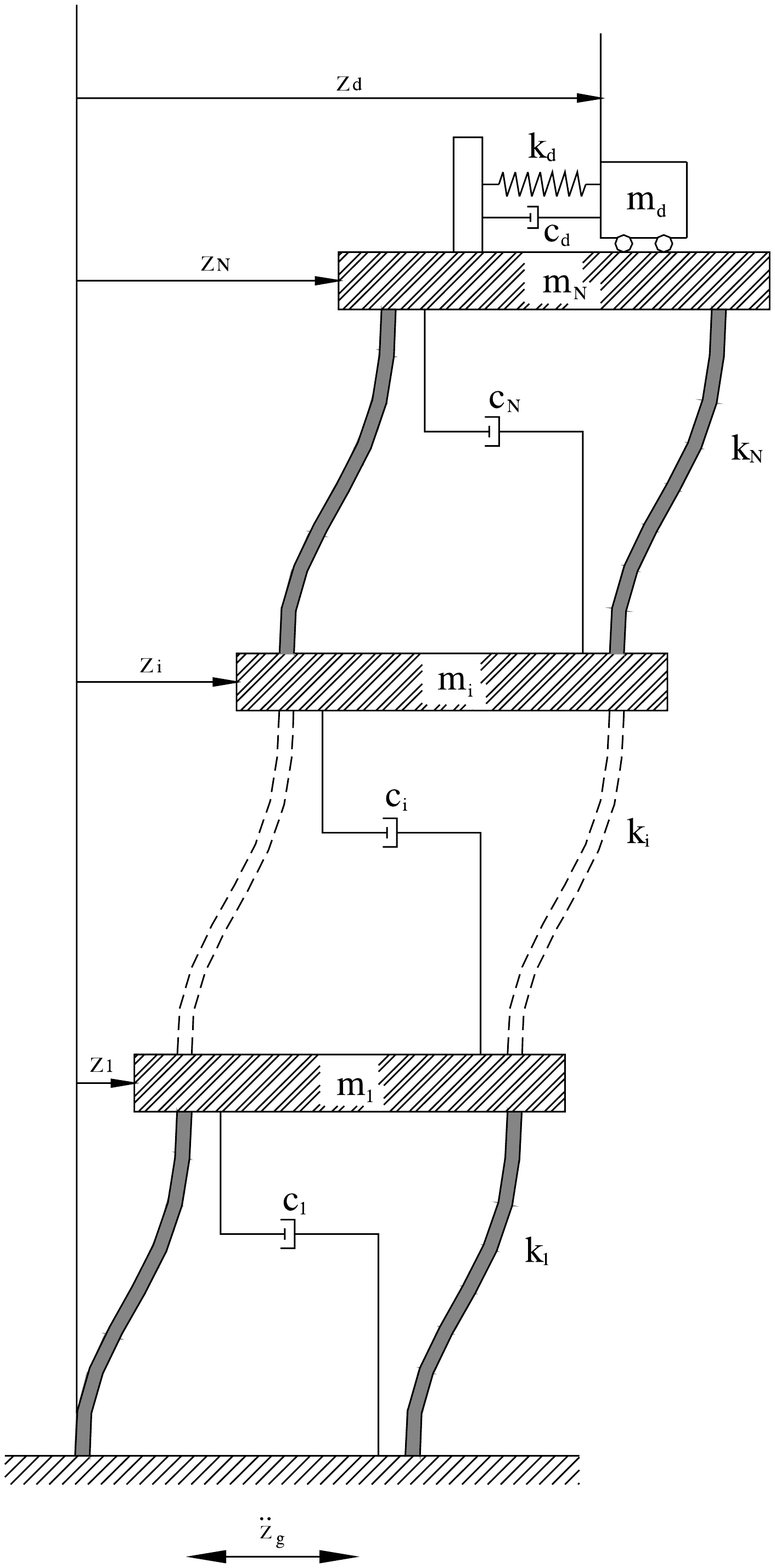}
    \caption{TDM building}
    \label{fig:tmd_building}
\end{figure}

Consider the $N$-story Multiple Degrees of Freedom (MDOF) linear building structure with a TMD installed at the top floor illustrated in \autoref{fig:tmd_building}. The equation of motion of the combined system subject to ground
acceleration can be written as 

\begin{equation}
    \mathcal{M} \ddot{\textbf{z}}(t)  + \mathcal{C} \dot{\textbf{z}}(t) + \mathcal{K} \textbf{z}(t)
    = - \textbf{m}~    \ddot{z}_g(t),
    \label{eq:motion}
\end{equation}

\noindent where $\textbf{z}$ is the $(N+1)$ dimensional response vector representing the displacements
relative to the ground

\begin{equation}
    \textbf{z}(t) = \{z_1,z_2,...,z_N,z_d \},
\end{equation}

\noindent $\ddot{z}_g$ is the ground acceleration, $\textbf{m}$ is the mass vector

\begin{equation}
    \textbf{m} = \{m_1,m_2,...,m_N,m_d\},
\end{equation}

\noindent and $\mathcal{M}$, $ \mathcal{C}$ and $\mathcal{K}$ are matrices corresponding to the
mass, viscous damping and the stiffness of the structure, respectively. These matrices can be written as

\begin{equation}
\small
    \mathcal{M} = \begin{bmatrix}
        m_1 & & \\
        & m_2 & \\
        & & \ddots \\
        & & & m_N \\
        & & & & m_d \\
    \end{bmatrix}
\end{equation}

\begin{equation}
\small
    \mathcal{C} = 
     \begin{bmatrix}
   (c_1+c_2) &  -c_2  &    &     &   \\
   -c_2 &  (c_2+c_3)  & -c_3  &     &   \\
    &  -c_3  &    &     &   \\
    &    &  \ddots  &  -c_N   &   \\
     &    &  -c_N  &  (c_N+c_d)   & -c_d  \\
     &    &    &  -c_d   & c_d  \\    
\end{bmatrix}
\end{equation}

\begin{equation}
\small
    \mathcal{K} = 
     \begin{bmatrix}
   (k_1+k_2) &  -k_2  &    &     &   \\
   -k_2 &  (k_2+k_3)  & -k_3  &     &   \\
    &  -k_3  &    &     &   \\
    &    &  \ddots  &  -k_N   &   \\
     &    &  -k_N  &  (k_N+k_d)   & -k_d  \\
     &    &    &  -k_d   & k_d  \\    
\end{bmatrix}
\end{equation} 

\noindent where $m_i$, $c_i$ and $k_i$ are respectively the mass, damping and stiffness of the $i$-th floor, while $m_d$, $c_d$ and $k_d$  are respectively the mass, damping and stiffness of the TMD.

We consider in this example that the structure is a ten-story ($N=10$) shear frame. The height of each floor is
3 meters, resulting in a total height (h) of 30 meters. Also, we model $m_i$, $c_i$ and $k_i$ ($i=1,...,N$) as Gamma random variables \citep{Ritto11} and group them into the random vector $\Xis$. Their mean values and Coefficient of Variation (C.o.V) are in
 \autoref{tab:statistical_parameters}. Notice here that the stochastic dimension of the problem is $n_x = 30$.

 \begin{table}[htbp]
    \small
    \centering
    \caption{Statistical properties of structural parameters}
    \renewcommand{\arraystretch}{1.5}
    \begin{tabular}{@{}cccc@{}}
        
        &       & Mean                & C.o.V. [\%]\\
        \toprule
        \multirow{2}{*}{\hfil Stifness [N/m]} & Story & $650.0 \times 10^6$ & 15 \\ \cline{2-4} & TMD & $k_d$               & 15 \\    \hline
        \multirow{2}{*}{\hfil Mass [kg] }     & Story & $360.0 \times 10^3$ & 05 \\ \cline{2-4} & TMD & $108.0 \times 10^3$ & 05 \\    \hline
        \multirow{2}{*}{\hfil Damping [Ns/m]} & Story & $6.20 \times 10^6$  & 25 \\ \cline{2-4} & TMD & $c_d$               & 25 \\
        \bottomrule
    \end{tabular}
    \label{tab:statistical_parameters}
 \end{table}
 
The earthquake excitation is modeled as a white noise signal with constant spectral density, $S_0$, filtered through the 
Kanai-Tajimi spectrum \citep{kanai,tajimi1960statistical}. The power spectral density function is given by:

\begin{equation}
s(\omega) = S_0 \left[ \frac{\omega_f^4+4\omega_f^2 \xi_f^2 \omega^2}
{(\omega^2-\omega_f^2)^2  + 4\omega_f^2 \xi_f^2 \omega^2} \right],
\end{equation}

\noindent where $\xi_f$ and $\omega_f$ are the ground damping and frequency, respectively. Their
values are adopted as $\xi_f = 0.6$, $\omega_f = 37.3$rad/s \citep{mohebbi2013designing}. The term
$S_0$ acts as a scaling factor and in this context represents the amplitude of the bedrock excitation
spectrum. Its value is adopted as $S_0= 1 \times 10^{-3} m^2/s^3$. This combination of parameters
corresponds to an earthquake with 0.38g peak ground acceleration on a medium firm soil
\citep{chen2005handbook}. The solution of the equation of motion is based on the Lyapunov equation. By solving
the  Lyapunov equation of the problem for the covariance matrix, it is possible to extract the
variance of the displacements and velocities from each degree of freedom. Therefore, it becomes
straightforward to  calculate the standard deviation of those quantities, which are in turn needed
for the reliability index computation. For a description of this procedure, the reader is referred to \citet{mantovani2017}.
 
The evaluation of the structural reliability ($\bar \beta$) leads to a time dependent reliability problem. Here, we employ the out-crossing rate approach for this purpose, which is detailed in Appendix \ref{sec:time_dependent_reliability}. For the calculation of the reliability index $\bar \beta$, the design life time ($t_D$) of the structure is
considered to be 50 years. Moreover, the rate of arrival ($\nu$) of earthquake events is of 1 every 10
years and  each event had the duration $t_E$ of 50 seconds. In this example, the top floor maximum displacement is chosen as failure criterion. Three different cases of failure barrier levels are analyzed: (a) $b = h/300$, (b) $b=h/400$, (c) $b=h/500$.

Thus, considering the design vector $\de = [k_d,c_d]$ and the stochastic parameters vector $\Xis$, the problem can be stated:

\begin{equation}
        \underset{\de \in S}{\min J(\de)} = - \bar \beta(\de),  \\ 
\end{equation}

\noindent in which the minus sign leads to the maximization of $\bar \beta$, which is given by \autoref{eq:fobj}. The design domain $S$, comprises the lower and upper bounds on the stiffness ($k_d$) and damping coefficient ($c_d$)
of the TMD. The bounds have values of $[0, 4000]$ kN/m and $[0,1000]$ kNs/m, respectively.

\begin{figure}[hbpt]
    \centering
    \includegraphics[width=\linewidth]{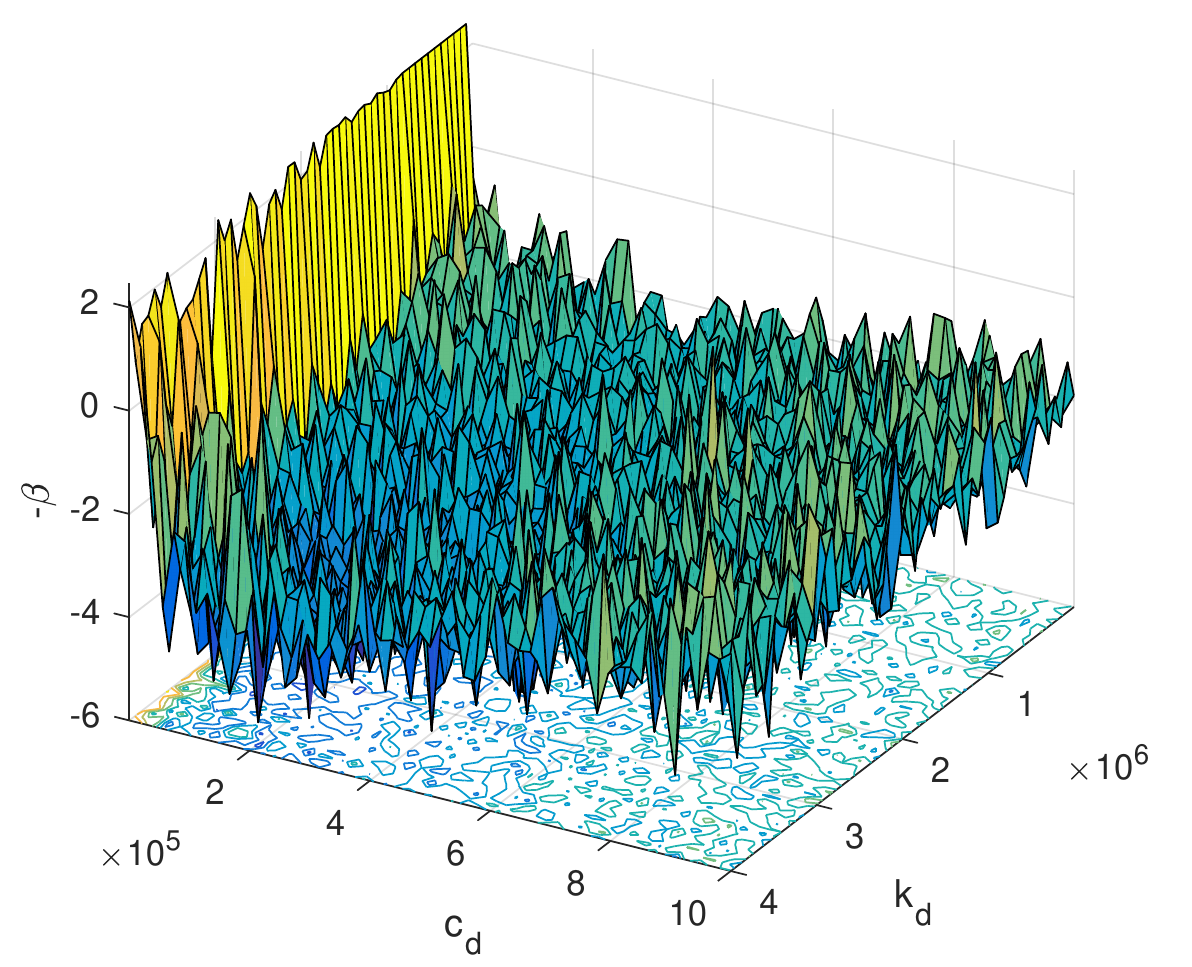}
    \caption{Noisy reliability surface over design variables range}
    \label{fig:surface_designvariables}
\end{figure}

In \autoref{fig:surface_designvariables} it can be seen how
irregular the surface becomes when considering the output without simulation, \ie the result of a
single evaluation for each input.

\autoref{tab:tmd_results} presents the optimization results using the proposed sEGO
approach with adaptive target variance. The stopping criterion of the algorithm is the maximum
number of calls of the FE code, \ie NFE$=1000$. The stiffness ($k_d$) and damping ($c_d$) are displayed along
with the best ($\bar \beta_{\text{max}}$) and average ($\bar \beta_{\text{ave}}$) results found over 25 independent runs of the algorithm. For comparison, the reliability index for the uncontrolled case is also displayed.

It can be seen that smaller barrier level results in smaller values of $\bar \beta$, thus increasing the failure probability. Notable increases in $\bar \beta$ were achieved. Taking for example the case (b),  the system reliability increases from 1.37 without TMD to 4.24 by using the TDM with the reported parameters. In terms of failure probability,
it means decreasing $P_f$ from $8.5 \times 10^{-2}$ to $1.1\times10^{-5}$. Looking at case (c), the structure that would certainly fail without TMD achieves a reliability index of 2.48 ($6.6\times10^{-3}$ failure probability) when using the optimized TMD parameters. Moreover, $\bar \beta_{\text{max}}$ results remained close to $\bar \beta_{\text{ave}}$. It indicates the robustness of the proposed approach, \ie it is able to obtain reasonable results in multiple runs despite the large number of random parameters and limited number of function evaluations.

\begin{table*}[htbp]
\centering
\caption{TMD optimization results}
\label{tab:tmd_results}
\begin{tabular}{@{}cccccc@{}}
\toprule
Barrier   & $k_d (MN/m)$ & $c_d (MNs/m)$ & $\beta_{max}$ & $\beta_{mean}$ & $\beta_{uncontrolled}$  \\ \midrule
(a) h/300 & 3.053        & 0.153         & 6.68          & 6.31           & 3.70            \\
(b) h/400 & 2.963        & 0.152         & 4.24          & 3.99           & 1.37            \\
(c) h/500 & 3.018        & 0.160         & 2.48          & 2.31           & fail            \\ \bottomrule
\end{tabular}
\end{table*}

\begin{figure}[ht]
    \centering
    \includegraphics[width=\linewidth]{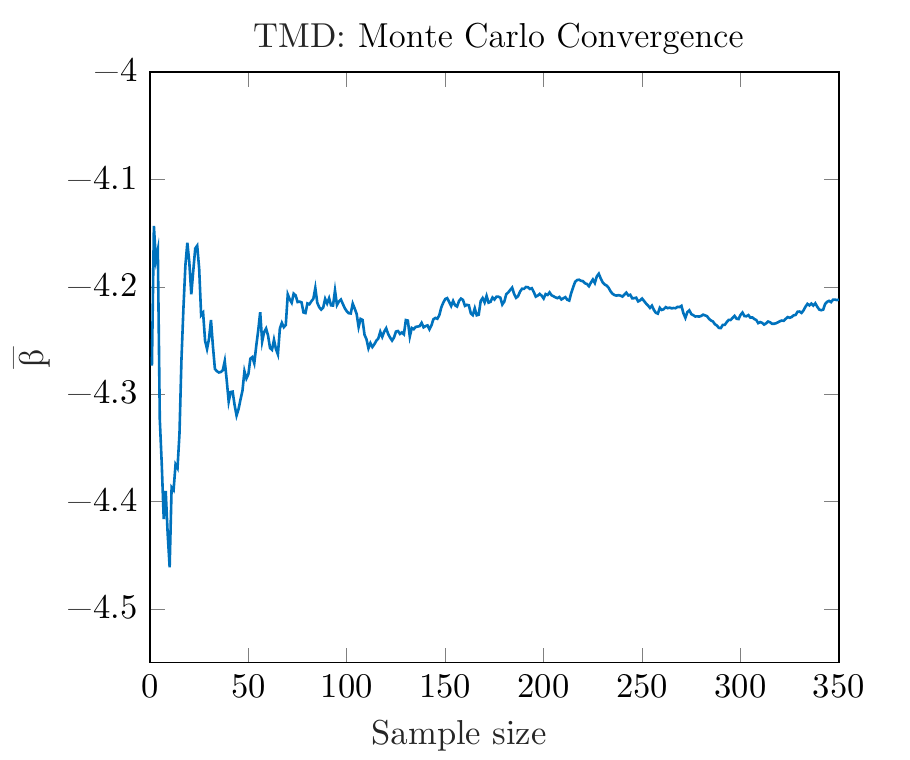}
    \caption{Mean convergence for the TMD reliability function}
    \label{fig:mcs_convergence}
\end{figure}

In \autoref{fig:mcs_convergence}, the Monte Carlo convergence for a single input is shown. The
average value of $-\bar \beta$ starts to converge around 150 simulations. Thus, it becomes clear that
applying the standard approach of using this fixed number to simulate every input would lead to a
higher computational cost. Simulating only seven different points would cause the maximum number of
evaluations of 1000 to be exceeded. However, using the proposed approach only points closer to
the optimum are simulated with higher $n_s$. Moreover, by employing the regression framework 
of SK using Monte Carlo variance estimates, one can approximate the underlying global behavior of 
the problem, avoiding further computational costs.

\begin{figure}[ht]
    \centering
    \includegraphics[width=\linewidth]{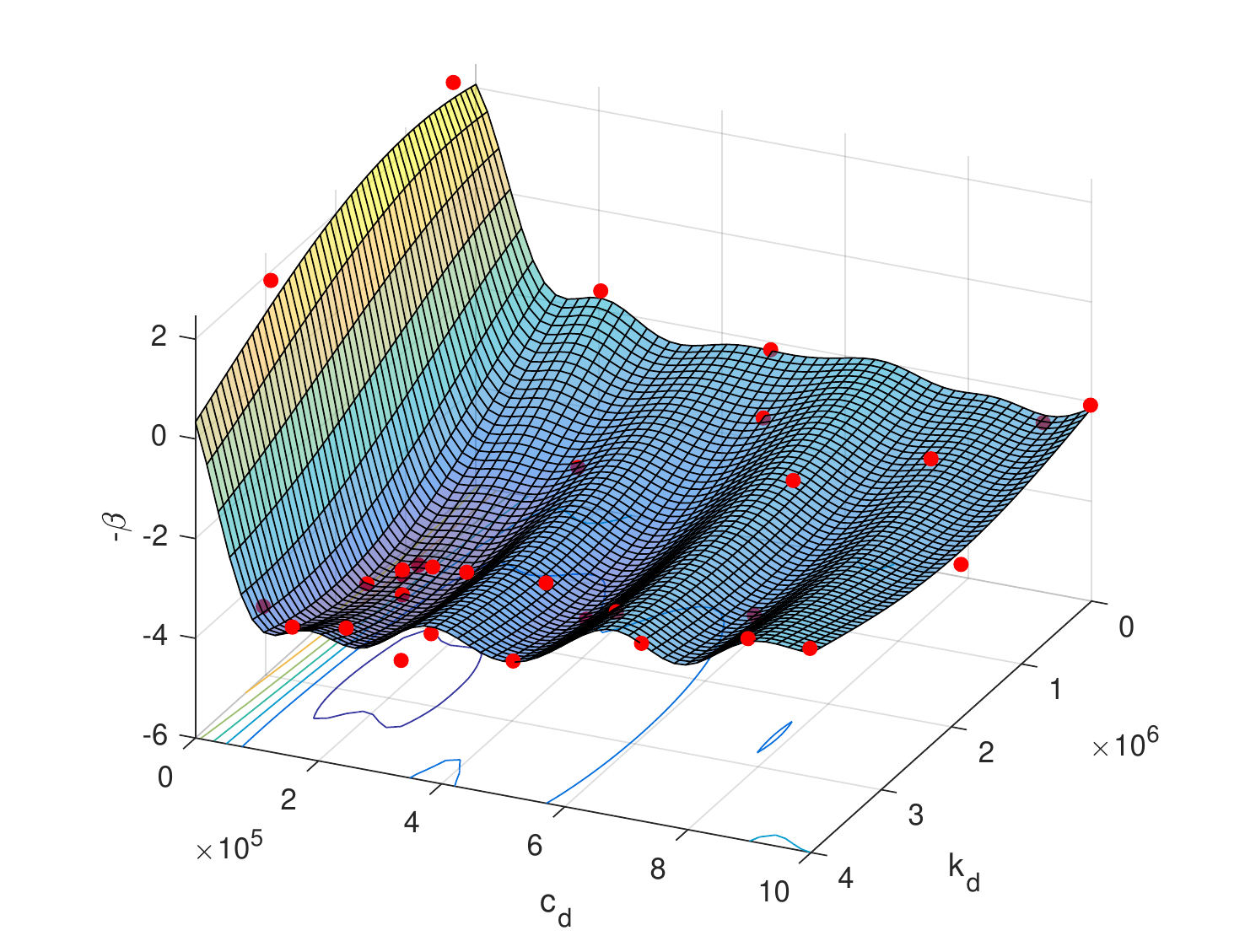}
    \caption{Surface generated by Stochastic Kriging with sampled points}
    \label{fig:surface_kriging}
\end{figure}

In \autoref{fig:surface_kriging}, the resulting surface generated by the proposed algorithm,
considering the case (b) barrier level is illustrated. Comparing it against \autoref{fig:surface_designvariables}
shows how the surface becomes smooth by the regression capabilities of SK. The red dots in
\autoref{fig:surface_kriging} represent the points that were sampled. Some are scattered over the
domain, which show the exploration part of the algorithm. Yet,
numerous points are concentrated close to minimum value of $-\bar \beta$, which exploited the region where the global optimum is located. Those are the points obtained
by the infill step, with the maximization of the AEI criterion and considering adaptive variance
target. This application possess the characteristics which the proposed algorithm is particularly
efficient at solving, \ie, problems with numerous random variables but a relatively small number of
design variables ($n\leq10$).

\subsection{Further comments}

The purpose of this section is offering general guidelines on parameter settings and discussing further application contexts. Overall most parameters remained constant 
throughout the analysis, where the biggest impact on the performance of the approach certainly comes from enforcing bounds on the target variance with $\overline{\sigma}^{2}_{target}$ and $\overline{\sigma}^{2}_{\min}$. The practical consideration is having an inexpensive upper bound and a limiting lower bound that assures that number of evaluations for each infill point can be completed and is compatible with the given computational budget. The actual values depend on characteristics of the problem being analyzed, where prior testing may be necessary to select reasonable parameters. The parameters of the exponential curve in \autoref{eq:target_adapt} remained constant across multiple dimensional problems and should not need major modifications. A side-length $r_{hc}$ of 10 \% of the range of each dimension seemed reasonable given that variables are normalized, performing well when the dimension size increased. The initial sampling plan, employed with $n_s=7n$ could be altered depending on the problem dimension and stochastic magnitude. For low error and low dimension, a reduced initial sampling plan would lead to a faster convergence, in reason of the higher probability of finding promising regions early.

Regarding further application, we mentioned in Section \ref{ModelSituation} that promising areas are robust design, Performance Based Design Optimization and Optimal Experimental Design. In the robust design of mechanical systems, MCI with fixed sample size has been largely employed to approximate the integral of Eq. (1), and metaheuristic algorithms have been widely applied to solve the resulting optimization problem. The coupling between fixed sample MCI and metaheuristic algorithms lead to a huge computational cost. In the case of expensive to evaluate functions, this is exactly what the global optimization approach proposed in this paper aims to avoid: the adaptive scheme avoids the full evaluation of the integral and the sEGO tends to find the global optimal solution requiring a fraction of the computational effort demanded by metaheuristic algorithms. Thus, robust design problems based on minimization of the mean response of the system would perfectly fit the capabilities of the algorithm proposed here.

Only a few papers in the literature addressed Performance Based Design Optimization problems, most likely for its huge computational cost (since it couples time dependent reliability analysis, nonlinear dynamics models and optimization). For example, \citet{beck2014optimal} simplified the problem to have only two random variables (applying a trapezoid rule for integration) and employed a local optimization algorithm. \citet{spence2014performance} built an equivalent static problem, but still had to solve a high dimension stochastic integral. The proposed sEGO approach could be employed in both situations: in the first, it would make it possible to consider more structural and loading parameters as random variables as well as to pursue a global optimization search, and in the second, it is likely to reduce the computational burden of the optimization process.

In the case of Optimal Experimental Design problems, we may also find some advantages of the proposed optimization method. The approach proposed by \citet{huan2013simulation} makes use of polynomial chaos surrogates, which may lose efficiency when the stochastic dimension of the problem increases (unless sparse approaches are employed). Thus, the adaptive approach would be preferred over polynomial chaos surrogates in situations with high stochastic dimension. On the other hand, \citet{Beck2017fast} employed a Laplace-based importance sampling method to alleviate the computational burden in the evaluation of the Optimal Experimental Design double integral. In this case, it is important to mention that this Laplace-based importance sampling method also provides the variance of the error in the integration. Hence, it can be employed (instead of MCI) together with the sEGO approach proposed here to further reduce the computational cost of the optimization problem. The same could be done with other efficient frameworks such as Multi Level Monte Carlo \citep{giles2008multilevel} and Multi index Monte Carlo \citep{haji2016multi}.


\section{Conclusion}
\label{sec:conclusions}

This paper presented an efficient sEGO approach for the minimization of functions that
depend on an integral. It was first supposed that this integral could be approximated by MCI, which also provided the variance of the error in the approximation. This information about the error was then included into the SK framework by setting a target variance in the MCI. The AEI
infill criterion was employed to guide the addition of new point in the metamodel. It was shown that if we set a large and constant target variance for MCI, it may stall the optimization process. On the other hand, if the target variance was relatively low, the computational cost might become prohibitive. It was then identified that there was room for optimization on the selection of the variance target. Hence, we then proposed an adaptive approach for variance target selection, which is the main contribution of this paper. 

In order to assess the effectiveness of the proposed method, numerous benchmark problems were
analyzed. The proposed method was first compared to the constant approach on deceptive stochastic benchmark functions from the
literature, with up to 10 design variables and up to 54 random variables. The proposed adaptive method obtained
better results in almost all analyzed cases. Moreover, it was observed that higher dimensions and higher noise levels led
to a higher difference in the results, favoring the proposed approach when compared to the constant target method. 

Then, the sEGO algorithms were compared to a multi-start global optimization method, the GBNM.  It was shown that the GBNM provided worse results than the sEGO approaches both in terms of efficiency, with a larger number of evaluations
required as well as consistency, with a higher 5 to 95 percentile range over multiple runs.

Finally, the design optimization under uncertainties of a TMD system was analyzed. It consisted of a ten-story shear frame subject to seismic excitation. The problem had thirty random variables and the objective was to
maximize the structural reliability by selecting the optimal stiffness and damping of the TMD system. The problem was successfully solved and the best results found were presented. The proposed adaptive approach was able to obtain optimal solutions efficiently and consistently for the three cases analyzed. Moreover, it supported its applicability to problems with high stochastic dimensions.

Overall, the proposed adaptive target sEGO method yielded convincing results for the minimization of
integrals. Furthermore, the stochastic dimensions of the integral to be minimized does not impose any limitation since we employed MCI for its approximation. By limitation, we mean that the method does not suffer from exponential increase in complexity when the number of stochastic dimensions is increased, also known as curse of dimensionality. Hence, the proposed sEGO method tends to be robust and efficient for the solution of problems that depend on the high dimensional integrals. On the other hand, one of the method limitation resides in the inherent size limitation of EGO coupled with the Kriging metamodel: the dimension $n$ of the optimization problem.


\begin{acknowledgements}
The authors acknowledge the financial support and thank the
Brazilian research funding agencies CNPq and CAPES.
\end{acknowledgements}

\bibliographystyle{plainnat}
\small
\bibliography{references}

\begin{appendices}

\section{Deterministic Kriging}
\label{dKriging}

Deterministic Kriging constructs a prediction model $\hat{y}$ based on the available information of the current sampling plan - $\boldsymbol{\Gamma}$ and $\textbf{y}$ - using Kriging \citep{sacks1989}. The basic idea behind Kriging is to construct a metamodel whose response at any point $\de$ \ is modeled as a
realization of a stationary stochastic process. Thus, at any point on the design domain, we have a Normal random variable with mean $\mu$ and variance $\sigma^2$. Considering an initial sampling plan $\boldsymbol{\Gamma}$, the covariance between any two input points $\de^{(i)}$ and $\de^{(j)}$ is:

\begin{equation}
\Cov\left[\de^{(i)},\de^{(j)}\right] = \sigma^2 \ \boldsymbol{\Psi}\left(\de^{
	(i)},\de^{
	(j)}\right),
\end{equation}

\noindent where $\boldsymbol{\Psi}$ is the correlation matrix, which has the form:

\begin{equation}
\boldsymbol{\Psi}\left(\de^{(i)},\de^{(j)}\right)  = \sum_{k=1}^{n} \exp \left(
-\theta_k \ \left|d^{(i)}_k - d^
{
	(j)}_k\right|^
{p_k} \right).
\label{eq:krigcorrelation}
\end{equation}

The unknown parameters $\theta_k$ and $p_k$ may be found by Maximum Likelihood Estimate (MLE) \citep{montgomery2010applied}, which then gives us the mean value - or average trend - and variance of the approximation: 

\begin{equation}
\hat{\mu} = \frac{\textbf{1}^T \boldsymbol{\Psi}^{-1} \textbf{y}}{\textbf{1}^T \boldsymbol{\Psi}^{-1} 
	\textbf{1}}
\label{eq:optimalmu}
\end{equation}

\noindent and

\begin{equation}
\widehat{\sigma^2} = \frac{(\textbf{y} - \textbf{1}\hat{\mu})^T \boldsymbol{\Psi}^{-1} (\textbf{y} -
	\textbf{1}\hat{\mu})}{n_s},
\label{eq:optimalvar}
\end{equation}

\noindent where $\textbf{1}$ is the identity matrix. With the estimated parameters, the Kriging prediction at a given point
$\de_u$ is:

\begin{equation}
\hat{y}(\de_u) =  \overbrace{\hat{\mu}}^{\text{Trend}} + \overbrace{\textbf{r}^T \boldsymbol{\Psi}^{-1} (\textbf{y} - \textbf{1} 
	\hat{\mu} )}^{\text{Model uncertainty}},
\label{eq:krigingpredictor}
\end{equation}

\noindent where \textbf{r} is the vector of correlations of $\de_u$ with the other $n_s$ Kriging sampled
points. The second term in the righ-hand-side of \autoref{eq:krigingpredictor} may be viewed as the model uncertainty since its value is inferred based on the function value of the points of the sampling plan.

One of the key benefits of kriging and other Gaussian process based models is the provision of an
estimated error in its predictions. The Mean Squared Error (MSE), derived by
\citet{sacks1989} using the standard stochastic process approach reads:

\begin{equation}
s^2(\de) = \widehat{\sigma^2} \left[ 1 - \textbf{r}^T \boldsymbol{\Psi}^{-1} \textbf{r} + \frac{(1 - \textbf{1}^T
	\boldsymbol{\Psi}^{-1} \textbf{r})^2}{\textbf{1}^T \boldsymbol{\Psi}^{-1} \textbf{1}} \right].
\label{eq:estimatederror}
\end{equation}

\autoref{eq:estimatederror} has the intuitive property that it is zero at already sampled points. In other words, the deterministic case of Kriging acts as a regression model which exactly interpolates the observed input/output data, \ie $\hat{y}(\de^{(i)}) = y^{(i)}$. 

Hence, deterministic Kriging assumes that the original model always provides an exact response, \ie $y^{(i)}=J(\de^{(i)})$. In other words, there is no error or variability when the original function $J$ is evaluated. However, it is not the case analyzed in this paper. Recall that we assume that Eq. (\ref{eq:problemabasico}) may not be analytically evaluated and we want to approximate it using MCI. Thus, for a given design vector $\mathbf{d}$, Eq. (\ref{eq:approxmean}) gives the approximation of the integral and Eq. (\ref{eq:pointestimate}) estimates the variance of the error in such an approximation. Consequently, when MCI is employed for the computation of $J$, the assumption made by deterministic Kriging no longer holds. Our goal is then to take this variability into account by constructing the metamodel using SK and including the information given by \autoref{eq:pointestimate} into the metamodel framework.

\section{Time-dependent reliability of oscillators}
\label{sec:time_dependent_reliability}

This section briefly presents the out-crossing approach for time dependent reliability problems. For a more detailed description, the reader is referred to \citet{Melchers}. The time-variant reliability problem for the random system response displacement can be formulated as follows. During a zero mean excitation event of specified duration $t_E$, the response of the oscillator should not exceed the specified limit - or barrier - $\pm b$. The barrier $b$ in this paper is the maxim displacement of top floor. If we have knowledge of system response statistics we may evaluate the out-crossing rate of the system, and consequently, its probability of failure. For a linear system excited by a zero mean Gaussian process, the response is Gaussian and the up-crossing rate can be
evaluated as

\begin{equation}
v_z^+(\mathbf{d},\mathbf{X}) = \frac{\sigma_{\dot{z}}(\mathbf{d},\mathbf{X})}{\sigma_{z}(\mathbf{d},\mathbf{X})}
\frac{1}{2 \pi} \exp\left(-\frac{b^2}{2 \left(\sigma_z(\mathbf{d},\mathbf{X})\right)^2} \right),
\label{eq:crossing}
\end{equation}

\noindent where $\sigma_z$ and $\sigma_{\dot{z}}$ are the standard deviation of the
displacement and of the velocity response, respectively. In Eq. (\ref{eq:crossing}), we make explicit the dependence of the crossing rate on the design variables $\mathbf{d}$ as well as the random  parameters $\mathbf{X}$ of the problem. Thus, considering a stationary excitation, the probability of a failure event $F$ of a given duration $ t_E $ may be computed as

\begin{align}
P(F|\mathbf{d},\mathbf{X},t_E)       & = 1 - \exp\left(-2  \int_0^{t_E} v_z^+(\mathbf{d},\mathbf{X}) dt \right), \\
& = 1 - \exp\left(-2 t_E v_z^+(\mathbf{d}, \mathbf{X}) \right).            \\
\end{align}


The structural loading from an earthquake, which is the application topic, is described by the arrival of an unknown number of events. Hence, the arrival of the events is modeled here as
a Poisson process. Thus, for a design life $t_D$ with a number of events $n_e$, we define the probability of failure $P_f$ as
\begin{align}
P_f(\mathbf{d},\mathbf{X}) & := P(F|\mathbf{d},\mathbf{X},t_D),          \label{eq:Pf}                                \\
& = \sum_{i=1}^\infty P(F|\mathbf{d},\mathbf{X}, t_E,n_e=i)~P(n_e=i|t_D),
\end{align}  
where                            
\begin{align}
P(F|\mathbf{d},\mathbf{X},t_E,n_e=i) & = 1 - (1- P(F|\mathbf{d},\mathbf{X},t_E))^i    \\
P(n_e=i|t_D)                 & = \frac{(\nu~t_D)^i \exp(-\nu~t_D)}{i!}, \label{eq:nu}
\end{align}

\noindent in which $\nu$ is the arrival rate of events. We may now obtain the reliability index $\beta$ from the $P_f$ given in Eq. (\ref{eq:Pf}) as 

\begin{align}
\beta(\mathbf{d},\mathbf{X})  & := - \Phi^{-1}(P_f(\mathbf{d},\mathbf{X})),       \label{eq:BetaX}
\end{align}

\noindent where $\Phi$ is the standard Gaussian cumulative distribution function.

%

%

%

Note that the reliability index in Eq. (\ref{eq:BetaX}) still depends on the random vector $\mathbf{X}$, characterized by its joint probability density function $f_{\mathbf{X}}$. Consequently, in order to compute the resulting structural reliability, we must then employ the Total Probability Theorem, leading to 

\begin{align}
\bar \beta(\mathbf{d}) & := \mathbb{E}_{\mathbf{X}} [\beta(\mathbf{d},\mathbf{X})]      \\
& = \int_{\mathbf{X}} \beta(\mathbf{d},\mathbf{x}) f_{\mathbf{X}}(\mathbf{x}) d\mathbf{x}  \label{eq:fobj} 
\end{align}

\noindent where $\mathbb{E}$ is the expected value operator and $\bar \beta(\mathbf{d})$ is the objective function to be maximized in the optimization process.

\section{Globalized Bounded Nelder–Mead (GBNM)}
\label{sec:gbnm}

The Globalized Bounded Nelder–Mead method employs a local search with a probabilistic restart, where the restart procedure uses an adaptive probability density function constructed using the memory of past local searches, as per \citet{luersen2004globalized}. The resulting optimization algorithm (coupling of the aforementioned restart and different local optimization algorithms) has been successfully applied to solve structural optimization problems (\citet{luersen2004globalized}, \citet{Ritto11}, \citet{TORII2011}) and it is described in the sequel.

In this paper, the local search employed was tha seme as in \citet{luersen2004globalized}, namely, Bounded Nelder-Mead. A starting point $\underline{\mathbf{s}}_0$ is chosen, then, a local search is performed and when a local minimum is found, the search is restarted. This restart procedure is described below.

The probability of having sampled a point $\underline{\mathbf{s}}$ is described by a Gaussian-Parzen-window approach \citet{duda2012pattern}:

\begin{equation}
f(\underline{\mathbf{s}})=\frac{1}{M}\sum\limits_{i=1}^M f_i(\underline{\mathbf{s}})\, ,\label{Eqprob1}
\end{equation}

\noindent where $M$ is the number of points $\underline{\mathbf{s}}_{(i)}$ already sampled. Such points come from the memory kept from the previous local searches, being, in the present version of the algorithm, all the starting points and local optima already found. $f_i(\underline{\mathbf{s}})$ is the Normal multivariate probability density function given by:

\begin{equation}
f_i(\underline{\mathbf{s}})=\frac{1}{(2\pi)^{n_p/2}\det{([\Sigma])}^{1/2}}\times\exp{\left(-\frac{1}{2}(\underline{\mathbf{s}}-\underline{\mathbf{s}}_{(i)})^T [\Sigma]^{-1}(\underline{\mathbf{s}}-\underline{\mathbf{s}}_{(i)})\right)}\,,
\end{equation}

\noindent where $n_p$ is the problem dimension and $[\Sigma]$ is the covariance matrix:

\begin{equation}
[\Sigma]=\left[
         \begin{array}{ccc}
           \sigma_1^2 &  &  \\
            & \ddots &  \\
            &  & {\sigma_{n_p}}^2 \\
         \end{array}
       \right]\,.
\end{equation}

The variances are estimated by the relation:

\begin{equation}
\sigma_j^2=\beta_o\left(s_j^{\text{max}}-s_j^{\text{min}}\right)^2
\end{equation}

\noindent where $\beta_o$  is a positive parameter that controls the length of the Gaussians, and $s_j^{\text{max}}$ and $s_j^{\text{min}}$ are the bounds of the $j^{th}$ variable ($j=1,2,..,n_p$). To keep the method simple, such variances are kept constant during the optimization process. At the end of each local search, $N$ points are randomly sampled $(\mathbf{s}_1,\mathbf{s}_2,\ldots,\mathbf{s}_N)$ and the one that minimizes Eq.(\ref{Eqprob1}) is selected as the initial point to restart the next local search. The stopping criterion of the global optimization of each subproblem is the maximum number of restarts ($nr_{max}$) defined \emph{a priori} by the user.

\end{appendices}

\clearpage

\end{document}